%% file: main_revised.tex
\documentclass[letterpaper, 11pt]{article} 
\usepackage{graphics,graphicx,xcolor,bm, enumerate,booktabs}
\usepackage{multicol} 
\usepackage{parskip}
\usepackage{amsmath,amssymb,mathtools,amsthm}
\usepackage{multirow,float}
\usepackage[utf8]{inputenc}
\usepackage{fancyhdr}
\usepackage[title]{appendix}
\usepackage{wasysym}
\usepackage{url}
\usepackage{arydshln}
\usepackage{ulem}
\usepackage{cancel}

\newtheorem{thm}{Theorem}[section]
\newtheorem{lemma}[thm]{Lemma}

\newtheorem{remark}{Remark}
\numberwithin{remark}{section}
\newtheorem{definition}{Definition}
\numberwithin{definition}{section}

\numberwithin{pro}{section}

\usepackage[font=footnotesize,labelfont=small]{caption}
\RequirePackage{geometry}
\newcommand{\lsp}{\vspace{3mm}}
\newcommand{\yzcmt}[1]{\textcolor{black}{#1}}
\newcommand{\agcmt}[1]{\textcolor{black}{#1}}

\newcommand{\mtx}[1]{\mathbf{#1}}
\newcommand{\vct}[1]{\mathbf{#1}}
\newcommand{\para}[1]{\textit{\textbf{#1.}}}

\title{ 
A fast direct solver for integral equations on locally refined boundary discretizations and its application to multiphase flow simulations  
}
\author{Yabin Zhang\thanks{Department of Mathematics, University of Michigan, Ann Arbor} , Adrianna Gillman\thanks{Department of Applied Mathematics, University of Colorado, Boulder}\, and Shravan Veerapaneni\textsuperscript{*}}

\date{\small The authors declared that they have no conflict of interest.}

\begin{document}

\maketitle

\begin{abstract}
In transient simulations of particulate Stokes flow, 
to accurately capture the interaction between the constituent particles and the confining wall,
the discretization of the wall often needs to be locally refined in the region approached by the particles. 
Consequently, standard fast direct solvers lose their efficiency since the linear system changes at each time step.
This manuscript presents a new computational approach that avoids this issue by pre-constructing a fast direct solver for the wall ahead of time, computing a low-rank factorization to capture the changes due to the refinement, and solving the problem on the refined discretization via a Woodbury formula. Numerical results illustrate the efficiency of the solver in accelerating particulate Stokes simulations.

\vspace{.5cm}

\end{abstract}
\section{Introduction}
\label{sec:intro}







A common computational task that arises in simulations of particulate Stokes flow is evaluating the hydrodynamic interaction of small moving geometries, such as drops, bacteria or biological cells, with large static structures, such as microfluidic chips, vascular walls, or channel walls. Boundary integral equation (BIE) methods, solved via iterative solvers accelerated by fast summation methods, are often used in practice for such systems as they avoid volume meshes as well as the cumbersome task of  volume re-meshing in transient simulations. In \cite{MARPLE2016_periodicstokes}, a fast direct solver was proposed which further reduces the cost of simulations by precomputing the compressed inverse of the BIE operator corresponding to the large static structures, which can be applied in linear time. This can be extremely useful in practice since most applications require a large number of time-steps to observe the physics of interest e.g., alignment of vesicles in a periodic channel \cite{ghigliotti2011vesicle}, pattern formation in suspensions of active particles \cite{lushi2014fluid, yan2020scalable} and cell sorting \cite{kabacaouglu2019sorting}.

However, when the suspended particles evolve in close proximity to the confining walls, the discretization of the walls must be locally refined to resolve the hydrodynamic interaction \cite{WU2020_adaptive}; this, in turn, makes direct solvers less attractive since the inverse operator needs to be re-evaluated continuously.  We present a fast algorithm that avoids re-building the inverse operator from scratch by precomputing an inverse operator corresponding to a reference mesh and rapidly updating it whenever the boundary discretization is locally refined (or coarsened). This work is an extension of Zhang-Gillman \cite{ZHANG2018, ZHANG2021}, where Laplace BIEs on locally-perturbed geometries were considered. The central idea is that the discretized BIE on the walls can be written as an extended version of the linear system for the original geometry and a fast direct solver on the original geometry can be reused to reduce the computational burden of solving the problems on the refined discretization. Since the conditioning of the discretized BIE for Stokes problem is at least the square of the Laplace BIE defined on the same geometry, special care is needed when using the Woodbury formula to apply the inverse of the extended system for numerical stability. 

\para{Related work}
At a high-level, fast direct solvers exploit
the fact that the off-diagonal blocks of the discretized system are low-rank.  
In the context of integral equations, some of them include the \textit{Hierarchically Block Separable (HBS)}  \cite{MARTINSSON20051,GILLMAN2012_china}, 
the
  \textit{Hierarchically Semi-Separable (HSS)}  \cite{Sheng2007_proceedings,CHANDRASEKARAN2004_withGu},
  the
   \textit{Hierarchical Interpolative Factorization (HIF)} \cite{HO2015_withYing_hif} and
  the $\mathcal{H}$ or $\mathcal{H}^2$- matrix methods \cite{HACKBUSCH1999_part1}.  
  The techniques developed in \cite{ZHANG2018, ZHANG2021} for the extended linear system (ELS), designed for problems with locally perturbed geometries, can be coupled with any of the above direct solver approaches. %
  \yzcmt{In this work, we employ a particular fast direct solver based on HBS matrix representation and inverse presented in \cite{GILLMAN2012_china}. 
  For the rest of the manuscript, when a HBS representation or inverse
  is built for a discretized boundary integral equation, it refers to
  the particular compression and inverse approximation given in \cite{GILLMAN2012_china}. } \agcmt{Other fast direct solvers can be used in place of the HBS solver and the results will be comparable.}

An alternative to using \agcmt{the} ELS is to update the hierarchical representation of the discretized integral operator directly.
Existing techniques in \cite{MINDEN2016,RYAN2020ARXIV} update the HIF of the system with \agcmt{a cost that is bounded above by the cost of } building a HIF of the 
perturbed or refined problem from scratch.
  For problems that do not
require a large number of discretization points, updating HIF directly is expected to be cheaper than 
building a new one from scratch.  
This idea is first investigated in \cite{MINDEN2016}, and a parallel implementation for Stokes BIEs on multiply-connected domains 
is presented in \cite{RYAN2020ARXIV}.
Being direct solvers, these techniques are
advantageous when a large number of solves are required for each new geometry.
Generalizing the idea to other standard fast direct solvers, such as those based on HBS or HSS matrix, 
requires knowledge of the particular compression techniques used in the chosen fast direct solver
and is non-trivial.

Several previous works employ fast direct solvers as preconditioners for the
linear systems that result from the discretization of integral equations and differential equations \cite{Darve_IFMM_pre,Hmat_pre,2010_HSS_pre,2005_beb,2003_beb}.  
Most of them build a low-accuracy direct solver for the linear
system and apply the forward operator via a fast matrix multiplication technique.
While convergence of the iterative solver is generally improved, it can be more dramatically
improved by the use of a more accurate direct solver as a preconditioner.
Section \ref{sec:preconditioner} explores the left preconditioner option
and how the accuracy of the direct solver impacts the quality of the preconditioner.

\para{Contributions} 
Motivated by the applications mentioned above, we apply the solution technique given in \cite{ZHANG2018,ZHANG2021} to Stokes flow problems defined on complex geometries,
some of which are adapted from real application geometry data.
\yzcmt{The linear system associated with the discretization of an integral equation for Stokes flow
has a physical nullspace corresponding to the pressure being unique up to a constant.
Fast direct solvers like HBS are sensitive to the existence of such nontrivial nullspace
due to the fact that matrices of smaller sizes are inverted in the hierarchical structure and singularity will immediately cause trouble.
The nullspace can be corrected via an analytic technique, but 
the resulting linear system can have high condition number due to the physics and/or
complexity of the geometry. In general, the linear system that
needs to be solved for Stokes problems have a condition number that is \textit{at least} squared that of the linear system  for a Laplace problem on the same geometry.
The high condition number of the system leads to similar condition number of the small
matrices inverted in the hierarchical structure of a fast direct solver,
resulting in loss of accuracy \agcmt{that is} not often seen in Laplace problems.
This is even more cumbersome when local refinement is added to the original discretization. 
The solution technique given in \cite{ZHANG2018,ZHANG2021} requires inverting a matrix 
whose conditioning may be worse than the original discretized BIE.   \agcmt{Since the 
condition number of the linear system for Stokes problems is often high (at least square that 
of Laplace), this technique without additional modifications to improve stability can be problematic.}}

\agcmt{The work in this manuscript improves the stability of the extended linear system solver from \cite{ZHANG2018,ZHANG2021}
by changing the technique used to create the low rank factorizations of the update matrix.  The updates to the previous versions
of the solver are inspired by the theory which specifies the conditions needed for the Woodbury formula to be stable.}
%
%
%
\agcmt{Even this updated solver cannot defeat a high condition number and how that impacts the accuracy of 
a direct solver.}
%
The work presented in this manuscript tackles this issue by
using the local refinement \agcmt{fast direct} solver as a preconditioner for the ELS.
When coupled with a
fast matrix multiplication technique for \agcmt{applying} the ELS, the resulting solution
technique converges in a constant number of iterations \textit{independent} of the number
of discretization points (as long as the geometric features are sufficiently resolved).


\para{Limitations} \agcmt{This manuscript only considers two dimensional problems even though}
the ideas introduced here generalize to higher dimensions.  \agcmt{Additional} work is needed in integrating with other computational machinery (e.g., quadratures) and carefully testing the efficiency of the overall solver.  
In dense suspension flows, the particle-wall near interactions happen over long length- and time-scales. Clearly, the solver developed here is not applicable to this setting since the wall geometry needs to be globally-refined, in which case the approach prescribed in \cite{MARPLE2016_periodicstokes} is better suited. Lastly, when the particles approach arbitrarily close to the walls, close evaluation schemes (e.g., \cite{barnett2015spectrally, WU2020_adaptive}) are required to improve the accuracy of interaction force computation. Although these schemes are not expected to change the computational efficiency, incorporating them and testing the solver is left to future work.  

%


\textit{\textbf{Outline.}}
The manuscript begins by reviewing boundary integral formulations for Stokes problems and
a technique for discretizing the resulting integral equations in Section \ref{sec:BIE}.  
Next the ELS for locally refined discretization and the corresponding
direct solver are presented in Section \ref{sec:fdsolver}.  The proposed preconditioner
for the ELS is presented in Section \ref{sec:preconditioner}.  
Next Section \ref{sec:numerics} illustrates the performance of the presented
solution techniques. Finally Section \ref{sec:conclusion} closes the manuscript
with a summary and concluding remarks.

\section{Boundary integral formulation}
\label{sec:BIE}
This manuscript considers integral equation techniques for solving both interior and exterior Stokes flow problems.  
The indirect integral equation formulation is employed, wherein, the solution can be cast as a convolution 
over the boundary $\Gamma$ of a kernel with an unknown boundary charge density.  For example, the 
velocity $\vct{u}$ can be represented by 
$$\vct{u} (\vct{x})=\int_\Gamma\mathcal{K}(\vct{x},\vct{y})\vct{\tau}(\vct{y}) ds_\vct{y} = (\mathcal{K}_\Gamma \vct{\tau})(\vct{x}),$$
where $\mathcal{K}$ denotes a kernel related to the fundamental solution of the Stokes equations and $\vct{\tau}$ denotes the unknown charge density.  
The kernel is chosen based on the problem under consideration.  One option is to represent the solution 
via the single layer integral operator denoted by $\vct{u}(\vct{x})  = (\mathcal{S}_\Gamma\tau)(\vct{x})$, where $\mathcal{S}$ denotes the Stokes single layer kernel (Stokeslet)
defined in its tensor components by

\begin{equation}
S_{ij}(\vct{x},\vct{y}) = \frac{1}{4\pi \mu} \left(\delta_{ij} \log\left(\frac{1}{r}\right) +\frac{r_ir_j}{r^2}\right), \quad i,j = 1,2, 
\label{eq:stokeslet}
\end{equation}
where  $\vct{r}:= \vct{x}-\vct{y}$,  $r = \|\vct{r}\|$ and $\delta_{ij}$ is the Kronecker delta.  

Another option is to use a double layer integral operator $\vct{u}(\vct{x})  = (\mathcal{D}_\Gamma\tau)$ to represent the velocity.  The tensor components of the double layer kernel $\mathcal{D}$ are
$$D_{ij}(\vct{x},\vct{y})= \frac{1}{\pi}\frac{r_ir_j}{r^2}\frac{\vct{r}\cdot \vct{n}_\vct{y}}{r^2}, \quad i,j = 1,2$$
where $\vct{n}_\vct{y}$ is the surface normal vector at the point $\vct{y}\in\Gamma$.  

Likewise, the pressure can be represented via an integral operator.  It should be chosen to 
match the representation of the velocity.  For example, if the velocity is represented with the single layer integral operator, then the 
pressure is given by 
$$p(\vct{x}) = \int_\Gamma \mathcal{Q}(\vct{x},\vct{y}) \vct{\tau}(\vct{y}) ds_\vct{y}$$
where 
$$Q_j(\vct{x},\vct{y}) = \frac{1}{2\pi}\frac{r_j}{r^2}, \quad j = 1,2$$
and $\vct{\tau}$ is the same boundary charge density as for the velocity. If the velocity is represented via the double layer integral operator, then the pressure is given by
$$p(\vct{x}) = \int_\Gamma \mathcal{P}(\vct{x},\vct{y}) \vct{\tau}(\vct{y})ds_\vct{y}$$
where $$P_j (\vct{x},\vct{y}) =\frac{\mu}{\pi}\left(-\frac{\vct{n}_{j,\vct{y}}}{r^2} + 
2\frac{r_j}{r^4}\vct{r}\cdot\vct{n}_\vct{y}\right), \quad j = 1,2$$
and $\vct{n}_{j,\vct{y}}$ denotes the $j^{\rm th}$ component of the surface 
normal vector $\vct{n}_\vct{y}$.


\subsection{Interior Stokes problem}
\label{sec:int}
Consider the  incompressible Stokes equation inside a geometry $\Omega_{\rm in}$ given by 

\begin{equation}\label{eq:model_bvp}
\begin{split}
-\mu \Delta \vct{u} (\vct{x})+ \nabla p (\vct{x}) &= \vct{0}, \quad   \ \mbox{ for } \vct{x}\in\Omega_{\rm in} \\
\nabla \cdot \vct{u}  (\vct{x})	      &= 0, \quad   \  \mbox{ for }\vct{x} \in\Omega_{\rm in} \\
 		   \vct{u} 	 (\vct{x})	      &=\vct{g} (\vct{x}),  \mbox{ for }\vct{x}\in\Gamma = \partial \Omega_{\rm in},
\end{split}
\end{equation}
where $\mu$ denotes the viscosity, $\vct{u}$ denotes the velocity, $\vct{g} (\vct{x})$ 
is a vector-valued function denoting the boundary data, and $p(\vct{x})$ is a scalar
valued function denoting the pressure. 
Figure \ref{fig:demo_figs}(a) gives a sample geometry.
The Dirichlet boundary data needs to satisfy the consistency condition
\begin{equation}
\int_{\Gamma}  \vct{g} (\vct{x})\cdot \vct{n}_\vct{x} \, ds_\vct{x} =0
\end{equation}
where $\vct{n}_\vct{x} $ denotes the outward pointing normal vector at $\vct{x}\in\Gamma$.

Representing the velocity via the double layer kernel
$$\vct{u}(\vct{x}) = (\mathcal{D}_\Gamma{\tau})(\vct{x})$$
results in having to solve the the following boundary integral equation
\begin{equation}\label{equ:interior_BIE}
-\frac{1}{2}\vct{\tau}(\vct{x}) +(\mathcal{D}{\vct{\tau}})(\vct{x})=g(\vct{x})
\end{equation}
for the unknown density $\vct{\tau}$ \cite{hsiao2008boundary}. 
Discretization of the BIE (\ref{equ:interior_BIE}) via the Nystr\"om 
method results in having to solve a dense linear system 
of the form 
\begin{equation}\label{equ:interior_BIE_discretized}
\left(-\frac{1}{2}\mtx{I} +\vct{D}\right)\bm{\tau}=\vct{g}xf
\end{equation}
where $\vct{D}$ denotes the matrix that results from the discretization of the 
double layer integral operator, $\vct{g}$ denotes a vector with entries given by the evaluation 
of $g(\vct{x})$ at the quadrature nodes, and the vector $\bm{\tau}$ 
denotes the vector of the unknown density values at the discretization points.

\begin{remark}
 The solution to (\ref{eq:model_bvp}) is unique up to a constant which 
results in the the linear system (\ref{equ:interior_BIE_discretized}) having a rank-1 nullspace.  
This nullspace can be 
corrected by adding the discretized integral operator $\mathcal{N}$
\begin{equation}\label{eq:null}
       (\mathcal{N}\tau)(\vct{x})=\vct{n}_\vct{x} \int_\Gamma \tau(\vct{y})\cdot \vct{n}_\vct{y}\, ds_\vct{y}
\end{equation}
to the discretized integral equation (\ref{equ:interior_BIE_discretized}).  
If a Nystr\"om discretization is used, the method in \cite{BIROS2004_stokes} can be used 
to discretize (\ref{eq:null}).  
Thus the linear system that needs to be solved for interior Stokes problems by 
using the double layer representation of the velocity is 
\begin{equation}\label{equ:interior_BIE_discretized2}
-\frac{1}{2} \bm{\tau} +(\vct{D} +\vct{N})\bm{\tau}=\vct{g}
\end{equation}
where $\vct{N}$ is the matrix that results from the discretization of (\ref{eq:null}).
 
\end{remark}

\subsection{Exterior Stokes problem}
\label{sec:ext}
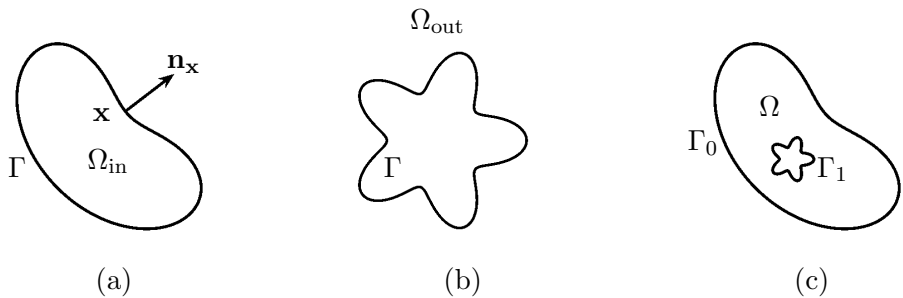
\begin{figure}[ht]
\begin{center}
\begin{tabular}{ccc}
\begin{picture}(120,100)(20,0)
\put(15,0){\includegraphics[height=30mm]{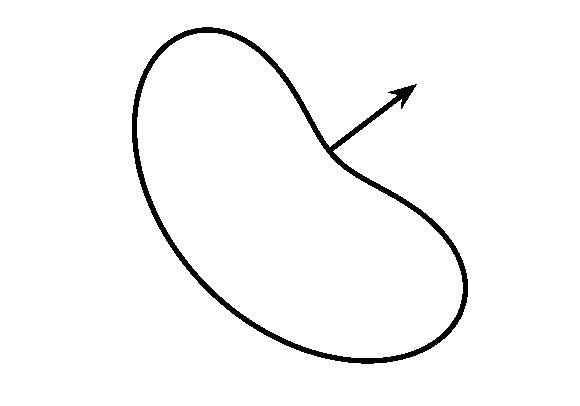}}
 \put(70,32){$\Omega_{\rm in}$}
\put(40,29){$\Gamma$}
\put(72,50){$\vct{x}$}
\put(100,70){$\vct{n}_{\vct{x}}$}
\end{picture}
&
\begin{picture}(120,100)(20,0)
\put(15,0){\includegraphics[height=30mm]{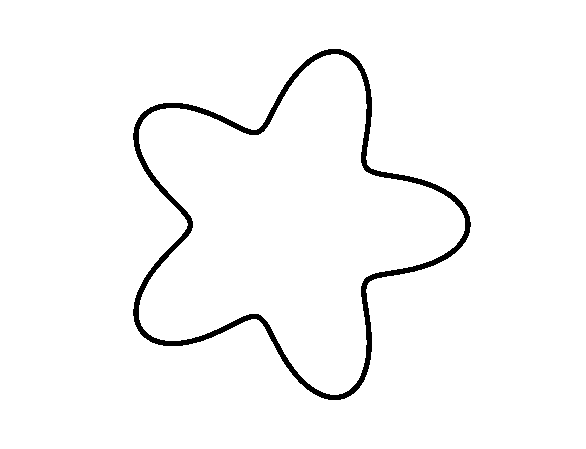}}
  \put(60,85){$\Omega_{\rm out}$}
\put(50,29){$\Gamma$}
\end{picture}
&
\begin{picture}(120,100)(20,0)
\put(15,0){\includegraphics[height=30mm]{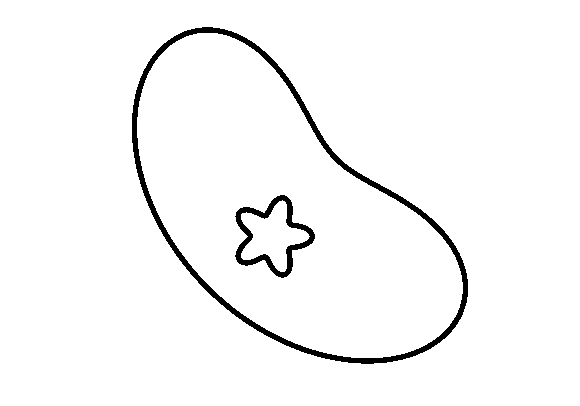}}
 \put(60,52){$\Omega$}
\put(33,39){$\Gamma_0$}
\put(82,29){$\Gamma_1$}
\end{picture}\\
(a) &(b) &(c)
\end{tabular}
\end{center}
\caption{ (a) A sample geometry for a purely interior BVP where the domain $\Omega_{\rm in}$ is the interior of
the boundary $\Gamma=\partial \Omega_{\rm in}$,
(b) a sample geometry for a purely exterior BVP where the domain $\Omega_{\rm out}$ is the exterior of
the boundary $\Gamma=\partial \Omega_{\rm out}$,
and (c) a sample geometry for an interior-exterior BVP where the domain $\Omega$ is the interior of the outer boundary $\Gamma_0$ 
but exterior of the inner boundary $\Gamma_1$.}
\label{fig:demo_figs}
\end{figure}

Exterior incompressible Stokes problems are also considered in this paper.  By an exterior 
problem, we mean that the velocity is sought in the domain $\Omega_{\rm out}$ defined 
as the plane minus the interior of a curve $\Gamma$ as shown in Figure \ref{fig:demo_figs}(b).
By using a combined field representation for the velocity
$$\vct{u} (\vct{x})= (\mathcal{D}_\Gamma\vct{\tau})(\vct{x}) +(\mathcal{S}_\Gamma\vct{\tau})(\vct{x}) = [\mathcal{(D+S)}_\Gamma \vct{\tau}](\vct{x}),$$
one is left with solving a second kind integral equation
\begin{equation}\label{equ:exterior_BIE}
\frac{1}{2}\vct{\tau}  +[\mathcal{(D+S)}_\Gamma \vct{\tau}]=\vct{g}.
\end{equation}
The linear system that results from discretizing this integral equation
is full-rank.

\begin{remark}
We also consider interior-exterior problems as shown in Figure \ref{fig:demo_figs}(c), where
the boundary $\Gamma$ is composed of an enclosing boundary curve $\Gamma_0$ and one or more holes with boundary $\Gamma_1$ inside the enclosed region.
The domain $\Omega$ is defined as the set that is interior to $\Gamma_0$ but exterior to $\Gamma_1$.
\end{remark}

\section{An extended linear system and direct solver for boundary value problems with locally refined discretization}\label{sec:fdsolver}
The efficient solver in this paper utilizes \agcmt{techniques previously developed}
 in \cite{ZHANG2018, ZHANG2021}, 
which are originally designed to handle BIEs defined on locally perturbed geometries. \agcmt{A geometry is said to be \textit{locally perturbed} if small parts of the boundary are modified from a previous BIE solve while the remainder of the boundary remains the same.}  We exploit the fact \agcmt{these} techniques can be applied to handle local discretization refinement.  \agcmt{For Stokes problems, the original fast solver techniques needed to be modified in order to handle
the higher condition number associated with these problems.  This section reviews the techniques from \cite{ZHANG2018,ZHANG2021}
and presents the new version needed for Stokes problems.}
Section \ref{sec:localperturbmodelproblem} begins by defining a 
 problem with locally refined discretization and introducing notation.  Section \ref{sec:localperturb_ext}
then presents the ELS and the efficient technique of solving 
that linear system using a solver built for the original discretization.
\yzcmt{
Section \ref{sec:low_rank_approx} introduces compression ideas for the blocks in the ELS that capture changes in discretization.
Finally,
Section \ref{sec:woodbury_stability} details the robustness of the solution technique and 
completes the algorithm.}

The fast direct solver presented in this section scales linearly with respect to the the number 
of points on the original discretization.  
\yzcmt{The solver can also scale linearly with 
respect to the number of points that are added in the refinement when a linear scaling inversion scheme is \agcmt{used
to invert} the discretized boundary integral operator on the refined part of the \agcmt{boundary}.
\agcmt{If the number 
of points added is not large (in general over a thousand), 
dense linear algebra is recommended for handling the refined region.  
This is because fast inversion algorithms such as HBS inversion \cite{MARTINSSON20051, GILLMAN2012_china}
tend to be slower than dense
linear algebra for small matrices.  }}

\subsection{Model problem with locally refined discretization}
\label{sec:localperturbmodelproblem}
\begin{figure}[ht]
\begin{center}
\begin{tabular}{cc}
\begin{picture}(150,130)(30,0)
\put(15,0){\includegraphics[height=45mm,trim={2cm 2cm 2cm 2cm}, clip]{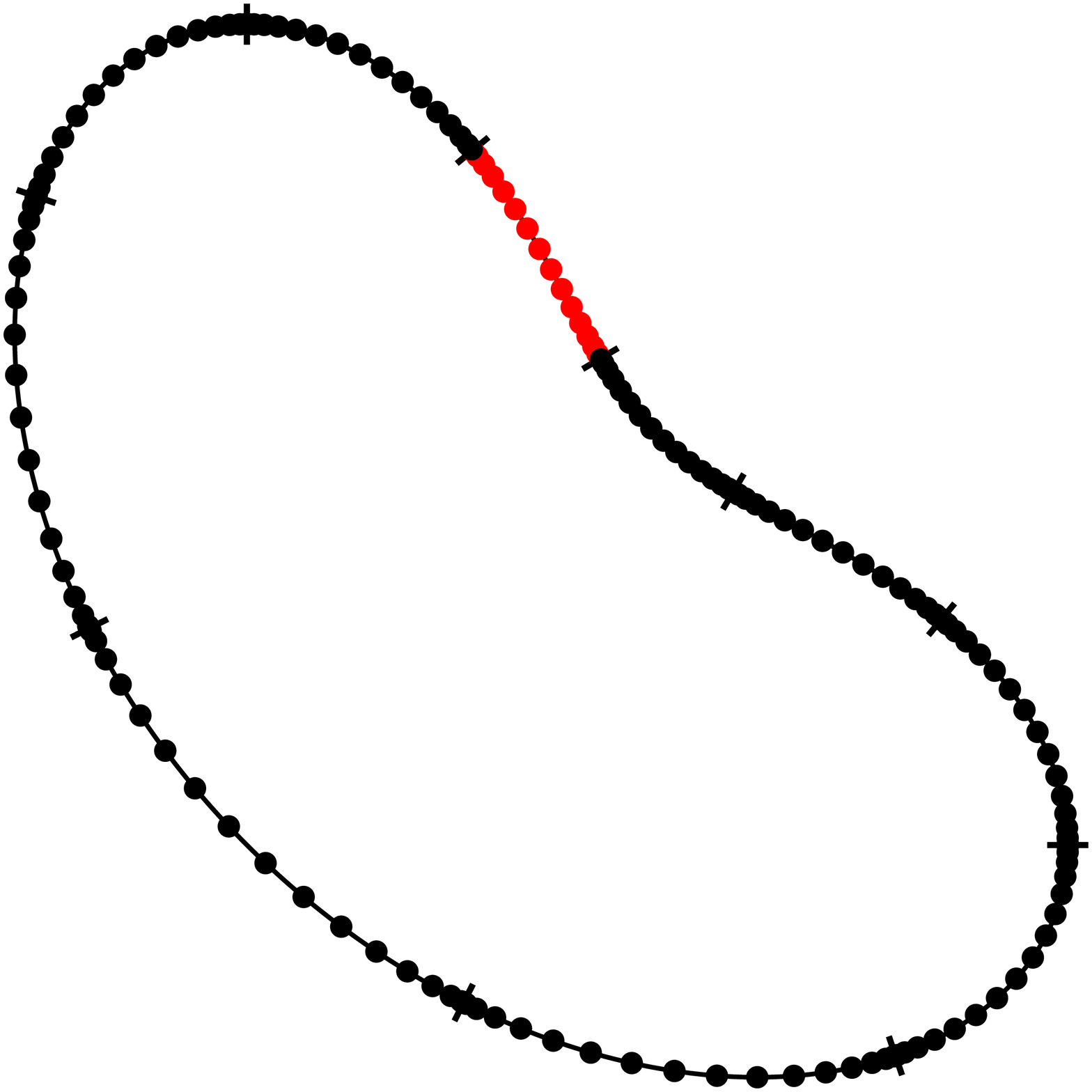}}
\put(50,29){$\Gamma_k$}
\put(110,100){$\Gamma_r$}
\end{picture}&
\begin{picture}(150,130)(10,0)
\put(15,0){\includegraphics[height=45mm,trim={2cm 2cm 2cm 2cm}, clip]{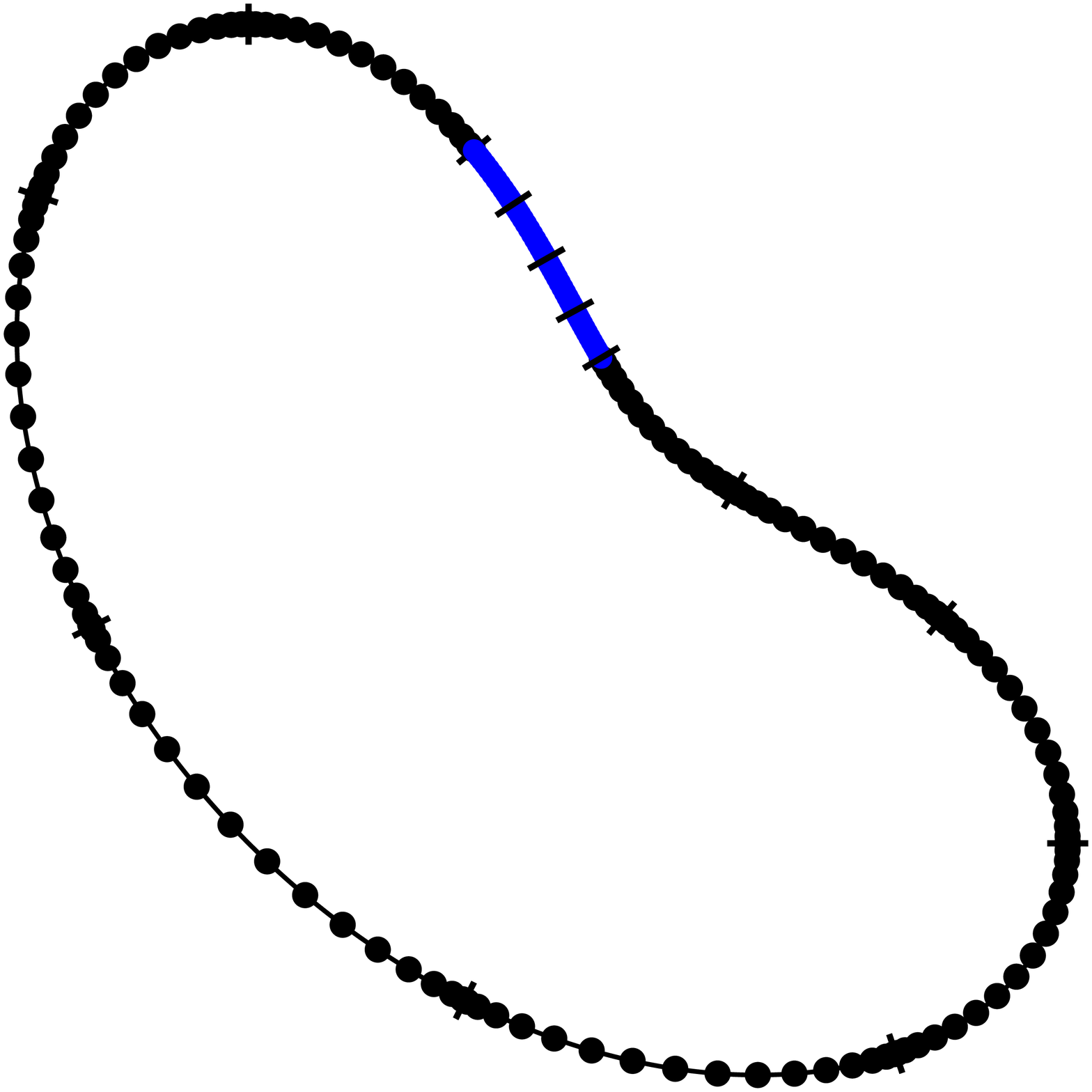}}
\put(50,29){$\Gamma_k$}
\put(110,100){$\Gamma_r$}
\end{picture}\\
(a)&(b)
\end{tabular}
\end{center}
\caption{ (a) The original Gaussian panel discretization of the geometry in Figure \ref{fig:demo_figs}(a). 
The discretization contains ten 16-point Gaussian panels uniformly distributed in parameterization space, and the panel in red is chosen to be refined.
(b) A locally refined discretization which \agcmt{replaced} the single red panel in Figure (a) with four blue panels. 
The part of the boundary curve that is refined is denoted by $\Gamma_r$ and the rest is denoted by $\Gamma_k$.
}
\label{fig:PerturbedBoundary}
\end{figure}

Consider the interior BVP defined by equation (\ref{eq:model_bvp}) on the geometry $\Omega_{\rm in}$ in Figure \ref{fig:demo_figs}(a).
As an example, let the boundary \agcmt{originally} be discretized with ten 16-point Gaussian panels.  \agcmt{Then one panel is chosen} to be refined into four panels.  See Figure \ref{fig:PerturbedBoundary}.
\agcmt{Let $\Gamma_r$ denote the part of the boundary that is refined and 
$\Gamma_k$ denote the part of the boundary where the discretization remains
 unchanged (``k" for ``kept").}
Figure \ref{fig:PerturbedBoundary}(a) and (b) \agcmt{illustrates} the pre- and post-refinement discretization respectively. 
The endpoints of the panels are also plotted.
For convenience, let $\mathcal{I}_k$, $\mathcal{I}_c$, and $\mathcal{I}_p$ denote the discretization points that are kept, deleted, and added for the refinement.
\agcmt{Thus}  $\mathcal{I}_o=\mathcal{I}_k\cup \mathcal{I}_c$ \agcmt{denotes} the collection of points \agcmt{in} the original discretization and \agcmt{$\mathcal{I}_n=\mathcal{I}_k\cup \mathcal{I}_p$ denotes the collection of discretization points on the boundary after refinement}.

The linear system (\ref{equ:interior_BIE_discretized2}) for the original and new discretization can be reordered in terms of the subscript notation.   Let $\mtx{A}$ denote the 
discretized integral equation (\ref{equ:interior_BIE_discretized2}) on the boundary; i.e., $\mtx{A}=-\frac{1}{2}\mtx{I}+\mtx{D}+\mtx{N}$.  
With the original discretization $\mathcal{I}_o$, 
the linear system can be ordered according to which points are added and deleted as follows
\begin{equation}\label{equ:originalLinearSystem}
\mtx{A}_{oo} \bm{\sigma}_{o}=
\begin{bmatrix}
\mtx{A}_{kk} & \mtx{A}_{kc}\\
\mtx{A}_{ck} & \mtx{A}_{cc}\\
\end{bmatrix}
\begin{pmatrix}
{\bm \sigma}_k\\
{\bm \sigma}_c
\end{pmatrix}
=
\begin{pmatrix}
\vct{g}_k\\
\vct{g}_c
\end{pmatrix}
=
\vct{g}_o.
\end{equation}
Likewise the linear system resulting from the refined discretization of the boundary integral equation can be ordered as follows
\begin{equation}\label{equ:perturbedlLinearSystem}
\mtx{A}_{nn} \vct{\bm{\tau}}_n =
\begin{bmatrix}
\mtx{A}_{kk} & \mtx{A}_{kp}\\
\mtx{A}_{pk} & \mtx{A}_{pp}\\
\end{bmatrix}
\begin{pmatrix}
\vct{\bm{\tau}}_k\\
\vct{\bm{\tau}}_p
\end{pmatrix}
=
\begin{pmatrix}
\vct{g}_k\\
\vct{g}_p
\end{pmatrix}
=
\vct{g}_n.
\end{equation}

\agcmt{In (\ref{equ:originalLinearSystem}) and (\ref{equ:perturbedlLinearSystem}), the subscript notation refers to the 
submatrices of $\mtx{A}$ on the respective geometries corresponding to different boundary interactions.  For example,
$\mtx{A}_{kk}$ denotes the submatrix of $\mtx{A}$ corresponding to the interaction of the points in $I_k$ with themselves ($A(I_k,I_k)$ in Matlab 
notation) and $\mtx{A}_{kp}$ denotes the submatrix of $\mtx{A}$ corresponding to the interaction of the points in $I_k$ with the points in $I_p$.}

\begin{remark}
 While the techniques in this section were presented for the interior problem (\ref{eq:model_bvp}), the techniques apply directly to exterior problems as well.
\end{remark}

\subsection{The extended linear system and direct solver}
\label{sec:localperturb_ext}
As an alternative to casting the problem solely on the ``new'' 
discretization, an ELS that is equivalent  to \agcmt{equation}
(\ref{equ:perturbedlLinearSystem}) can be considered.  In this 
paper, we use the ELS from \cite{ZHANG2021}.  
The ELS 
takes the form 
\begin{equation}\label{eq:extendedlLinearSystem}
\mtx{A}_{\rm ext} \bm{\tau}_{\rm ext} =
\begin{bmatrix}
\mtx{A}_{kk} & \mtx{0} & \mtx{A}_{kp}\\
\mtx{A}_{ck} & \mtx{A}_{cc}  & \mtx{0}\\
\mtx{A}_{pk} & \mtx{0} & \mtx{A}_{pp}\\
\end{bmatrix}
\begin{pmatrix}
&{\bm{\tau}}_k\\
&{\bm{\tau}}^{\rm dum}_c\\
&{\bm{\tau}}_p
\end{pmatrix}
=
\begin{pmatrix}
\vct{g}_k\\
\vct{0}\\
\vct{g}_p
\end{pmatrix}
=\vct{g}_{\rm ext}
\end{equation}
where ${\bm{\tau}}_k$ and ${\bm{\tau}}_p$ are the unknown boundary 
densities evaluated at points in $\mathcal{I}_n$
and ${\bm{\tau}}^{\rm dum}_c$ is 
a dummy boundary density \agcmt{at the points in $\mathcal{I}_c$} that is not used to evaluate the solution 
in the domain.  This linear system can be written 
as $\mtx{A}_{\rm ext}= \tilde{\mtx{A}}+ \mtx{Q}$ where
 $$
\tilde{\mtx{A}}=\begin{bmatrix}
\mtx{A}_{oo} & \mtx{0} \\
 \mtx{0} & \mtx{A}_{pp}\\
\end{bmatrix}
\;\mbox{ and }
\mtx{Q}=
\begin{bmatrix}
 \mtx{0}  & -\mtx{A}_{kc} & \mtx{A}_{kp}\\
 \mtx{0} &  \mtx{0} & \mtx{0}\\
\mtx{A}_{pk} & \mtx{0} & \mtx{0}\\
\end{bmatrix}.
$$
The matrix $\tilde{\mtx{A}}$ is full rank and block-diagonal with the first block 
equal to the operator for the original discretization.  Thus if the inverse of $\mtx{A}_{oo}$ 
has been precomputed (directly or via a fast direct solver), the cost of inverting 
$\tilde{\mtx{A}}$ is the cost of the inverting $\mtx{A}_{pp}$ which is small in
the problems under consideration.  The update matrix $\mtx{Q}$ is a block sparse 
matrix consisting of only three non-zero sub-blocks.  Since these non-zero blocks 
of $\mtx{Q}$ correspond to non-self interactions, they are low rank; i.e., $\mtx{Q}$ 
is low rank.  Let $\mtx{Q} = \mtx{LR}$ denote the low rank factorization of $\mtx{Q}$.  

The advantage of writing the linear system in the extended form (\ref{eq:extendedlLinearSystem}) 
and writing it as the sum of a block diagonal matrix with a low rank matrix is that the 
 inverse can be approximated via a Woodbury formula
\begin{equation}\label{eq:wood}
\bm{\tau}_{\rm ext}=
\left( \tilde{\mtx{A}} +\mtx{Q} \right)^{-1}\vct{g}_{\rm ext}\approx
\left( \tilde{\mtx{A}} +\mtx{LR} \right)^{-1}\vct{g}_{\rm ext}
\approx
\tilde{\mtx{A}}^{-1}\vct{g}_{\rm ext}
-\tilde{\mtx{A}}^{-1}\mtx{L}
  \left( \mtx{I} + \mtx{R}\tilde{\mtx{A}}^{-1}\mtx{L}\right )^{-1}
  \mtx{R}\tilde{\mtx{A}}^{-1}\vct{g}_{\rm ext}.
\end{equation}
This inverse can be applied rapidly to vectors by exploiting the block structure
of the matrices.  The only matrix that needs to be inverted in the application of (\ref{eq:wood}) is  $ \left( \mtx{I} + \mtx{R}\tilde{\mtx{A}}^{-1}\mtx{L}\right )$. This matrix is of size $k=k_{kc}+ k_{kp} +k_{pk}$ where 
$k_{kc}$, $k_{kp}$, and $k_{pk}$ denote the $\epsilon-$ranks of the low-rank approximations of $\mtx{A}_{kc}$, $\mtx{A}_{kp}$, and $\mtx{A}_{pk}$ given tolerance $\epsilon$, respectively.  Typically, $k$ is small and thus the matrix can 
be inverted via dense linear algebra for little computational cost. 
Algorithm  \ref{alg:fdsolver} summarizes the technique for rapidly applying the inverse 
of $\mtx{A}_{\rm ext}$ provided a fast direct solver for $\mtx{A}_{oo}$ has already been 
computed.  The algorithm is designed so that it can be used with any fast direct 
solver including the HBS\cite{GILLMAN2012_china}, HSS\cite{Sheng2007_proceedings,CHANDRASEKARAN2004_withGu}, HIF\cite{HO2015_withYing_hif}, and $\mathcal{H}$ or $\mathcal{H}^2$- matrix methods\cite{HACKBUSCH1999_part1}.
\yzcmt{Section \ref{sec:low_rank_approx} presents fast techniques for \agcmt{creating the low rank factorizations of} the blocks in $\mtx{Q}$ and section \ref{sec:woodbury_stability} discusses the stability for using the Woodbury formula and necessary improvements for the \agcmt{low rank factorization} of $\mtx{Q}$.
\begin{remark}
The  factorization technique for $\mtx{Q}$ (step 1 in Algorithm 1) to be discussed in section  \ref{sec:low_rank_approx} and \ref{sec:woodbury_stability} scales linearly with respect to $N_k$, $N_c$ and $N_p$.
Thus, if a fast direct solver is constructed for $\mtx{A}_{pp}^{-1}$, 
then all steps in pre-compuation and solve of Algorithm 1 scale linearly with respect to $N_k$, $N_c$ and $N_p$.
Otherwise, Algorithm 1 scales linearly with respect to $N_k$ and $N_c$ but cubically with respect to $N_p$ due to the dense linear algebra calculations for $\mtx{A}_{pp}^{-1}$. Table \ref{tab:scaling_compare} lists the cost scaling of Algorithm 1 and the fast direct solver in \cite{ZHANG2018}. 
More details on the step-by-step cost analysis is given in \cite{ZHANG2018}.
The scaling for the fast direct solver given in \cite{ZHANG2021} is the same as  Algorithm 1. 
\end{remark}
}

\begin{table}
    \centering
    \begin{tabular}{|l|c|c|}
    \hline
     Method    &   Pre-computation & Solve \\
         \hline
    \cite{ZHANG2018} with dense linear algebra for $\mtx{A}_{pp}$    & $O\left(N_k+N_c^3+N_p^3\right)$& $O\left(N_k+N_c^2+N_p^2\right)$\\
    \cite{ZHANG2018} with fast direct solver for $\mtx{A}_{pp}$    & 
    $O\left(N_k+N_c^3+N_p\right)$& $O\left(N_k+N_c^2+N_p\right)$\\
    \hline
    Algorithm 1 with dense linear algebra for $\mtx{A}_{pp}$    &
    $O\left(N_k+N_c+N_p^3\right)$&$O\left(N_k+N_c+N_p^2\right)$\\
    Algorithm 1 with fast direct solver for $\mtx{A}_{pp}$    &
    $O\left(N_k+N_c+N_p\right)$& $O\left(N_k+N_c+N_p\right)$\\
    \hline
    \end{tabular}
    \caption{Cost scaling for the fast direct solver in \cite{ZHANG2018} and the proposed solver in Algorithm 1. The fast direct solver in \cite{ZHANG2021} has the same scaling as Algorithm 1.}
    \label{tab:scaling_compare}
\end{table}


 

\begin{figure}
\begin{center}
\fbox{
\begin{minipage}{.9\textwidth}
\begin{center}
\textsc{Algorithm  \ref{alg:fdsolver}}: {\agcmt{Applying the} fast direct solver \agcmt{for} the locally refined problem }
\end{center}

\lsp

\textit{Given a fast direct solver for the original discretization $\mtx{A}_{oo}^{-1}$, and the right-hand-side vector defined for the refined discretization $\vct{g}_n=\begin{pmatrix}\vct{g}_k\\
\vct{g}_p \end{pmatrix}$,  
 this algorithm 
determines the solution to the refined problem (\ref{equ:perturbedlLinearSystem}) by obtaining the solution to the equivalent ELS via a Woodbury formula (\ref{eq:wood}).
}

\lsp

 \begin{tabbing}
  \hspace{5mm} \= \hspace{5mm} \= \hspace{5mm} \= \hspace{60mm} \= \kill
  
  \textit{Pre-computation:}\\
 \> \textbf{Step 1:} \yzcmt{Factorize the update matrix $\mtx{Q}\approx \mtx{LR}$ via the method in Section \ref{sec:low_rank_approx} and Section \ref{sec:woodbury_stability}.}\\
  \> \textbf{Step 2:} (invert $\mtx{A}_{pp}$)\\
  \>\textbf{if} $N_p$ is small,\\
  \> \> \yzcmt{Evaluate and invert $\mtx{A}_{pp}^{-1}$ via dense linear algebra.}\\
  \> \textbf{else},\\
  \> \>  \yzcmt{Build \agcmt{an approximate inverse of $\mtx{A}_{pp}$ via a fast direct solver such as HBS}.}\\
  \> \textbf{end if}\\
  \>  \textbf{Step 3:} \yzcmt{Apply the applying scheme for $\mtx{A}_{oo}^{-1}$ and  $ \mtx{A}_{pp}^{-1}$ to evaluate $\mtx{X}=\tilde{\mtx{A}}^{-1}\mtx{L}$.}\\
  \>  \textbf{Step 4:} \yzcmt{Evaluate and invert the Woodbury operator $ \left( \mtx{I} + \mtx{R}\mtx{X}\right )$  via dense linear algebra.}\\
\>\\
\>\\  
  \textit{Solve:}\\
  
   \>  \textbf{Step 1:} Evaluate $\tilde{\mtx{A}}^{-1}\bm{g}_{\rm ext}
  =  \begin{pmatrix}
\mtx{A}_{oo}^{-1}\begin{pmatrix}\vct{g}_k\\
\vct{0}\end{pmatrix}\\
\mtx{A}_{pp}^{-1}
\vct{g}_p
\end{pmatrix}$ utilizing the \agcmt{fast matrix vector applies provided}\\
\> \hspace{1.45cm}  \agcmt{by the direct solver(s)}.\\
 \>  \textbf{Step 2:} Evaluate $\bm{\tau}_{\rm ext}$ via the Woodbury formula (\ref{eq:wood}).\\
 \end{tabbing}
\end{minipage}}
\end{center}
\end{figure}


\subsection{Efficient construction of the low-rank factors in $\mtx{Q}\approx \mtx{LR}$}
 \label{sec:low_rank_approx}
 \yzcmt{
 There are two steps in the proposed technique for creating the low-rank approximation of $\mtx{Q}$.
This section introduces the first step
which creates low-rank factorizations for the non-zero blocks in $\mtx{Q}$. 
The second step, a recompression step
which is necessary for avoiding conditioning issues
associated with using the Woodbury formula to apply the inverse of the 
ELS in Equation (\ref{eq:wood}), is the delayed to the next section
after a brief review on the numerical stability of the Woodbury formula.
}

The low rank factorization of the update matrix $\mtx{Q}$ is done in a block format.
In other words, low rank factorizations are constructed for each of the non-zero subblocks
of $\mtx{Q}$; 

\begin{center}
\begin{tabular}{ccclcccclc}
$\mtx{A}_{kc}$ & $\approx$ &$ \mtx{L}_{kc}$ &$\mtx{R}_{kc}$, &       &$\mtx{A}_{kp} $&  $\approx$ & $\mtx{L}_{kp}$ &$\mtx{R}_{kp},$ \qquad $\mbox{ and }$&\\ 
{\footnotesize$ 2N_k \times 2N_c $}& & {\footnotesize$2N_k\times k_{kc}$} &  {\footnotesize$k_{kc}\times 2N_c$} &   & {\footnotesize$2N_k \times 2N_p$ }& & {\footnotesize$2N_k\times k_{kp}$} & {\footnotesize$k_{kp}\times 2N_p$}&\\
\end{tabular}
\end{center}
\begin{equation}\label{equ:block_lr}
\begin{tabular}{cccl}
$\mtx{A}_{pk} $&  $\approx$ & $\mtx{L}_{pk}$ &$\mtx{R}_{pk}. $\\
{\footnotesize$2N_p \times 2N_k$} & & {\footnotesize$2N_p\times k_{pk} $}& {\footnotesize$ k_{pk}\times 2N_k$}\\
\end{tabular}
\end{equation}
Here $N_k$, $N_c,$ and $N_p$ are the number of discretization points in $\mathcal{I}_k, \mathcal{I}_c,$ and $\mathcal{I}_p$ respectively.

Thus the low-rank factorization of $\mtx{Q}$ can be expressed as
\begin{equation}
\label{equ:QeqLR}
\begin{tabular}{cccl}
$\mtx{Q}$ &$\approx$ &$\mtx{L}_1$ & $\mtx{R}_1$\\
{\footnotesize$2N_{\rm ext}\times 2N_{\rm ext}$}& & {\yzcmt{\footnotesize$2N_{\rm ext} \times k_1$}} & {\yzcmt{\footnotesize$k_1 \times 2N_{\rm ext}$}}
\end{tabular}
\end{equation}
where
\begin{equation*}
\mathbf{L}_1=
\begin{bmatrix}
\mtx{0} & -\mtx{L}_{kc} & \mtx{L}_{kp}\\
\mtx{0} &\mtx{0} & \mtx{0}\\
\mtx{L}_{pk} &\mtx{0}&\mtx{0}\\
\end{bmatrix}  
\mbox{, }
\mathbf{R}_1=
\begin{bmatrix}
\mtx{R}_{pk} & \mtx{0} &\mtx{0}\\
\mtx{0} & \mtx{R}_{kc} &\mtx{0}\\
\mtx{0}  & \mtx{0} & \mtx{R}_{kp}
\end{bmatrix},
\end{equation*}
 $k_1=k_{pk}+k_{kc}+k_{kp}$ and $N_{\rm ext}=N_k+N_c+N_p$. 
 \yzcmt{Note the subscript notation in $\mtx{L}_1, \mtx{R}_1$ and $k_1$
 are intended because we reserve the notation   $\mtx{L}$, $\mtx{R}$ and $k$
 for the final low-rank factorization of $\mtx{Q}$ obtained \agcmt{from} the recompression \agcmt{which is presented} in section \ref{sec:woodbury_stability}. }
 
\begin{figure}[ht]
\begin{center}
\begin{tabular}{ccc}
 \begin{picture}(150,130)(30,0)
\put(15,0){\includegraphics[height=45mm,trim={2cm 2cm 2cm 2cm}, clip]{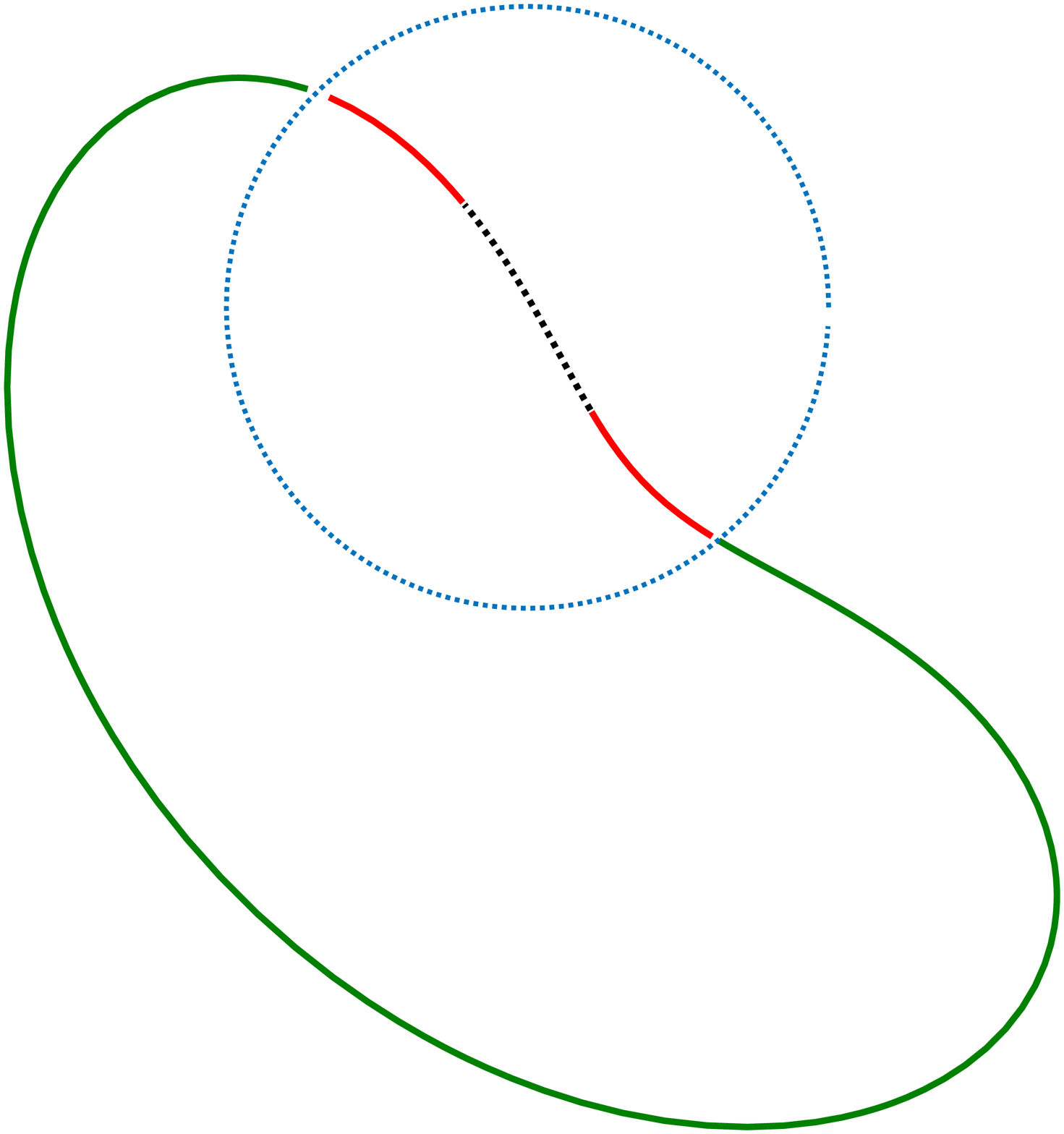}}
\put(110,95){$\Gamma_r$}
\end{picture}&
\begin{picture}(50,130)(0,0)
\put(0,30){\includegraphics[height=25mm,trim={0cm 0cm 0cm 0cm}, clip]{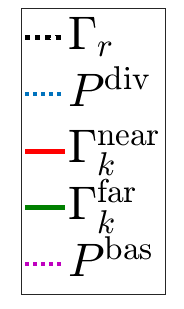}}
\end{picture}
&
\begin{picture}(150,130)(30,0)
\put(15,0){\includegraphics[height=45mm,trim={2cm 2cm 2cm 2cm}, clip]{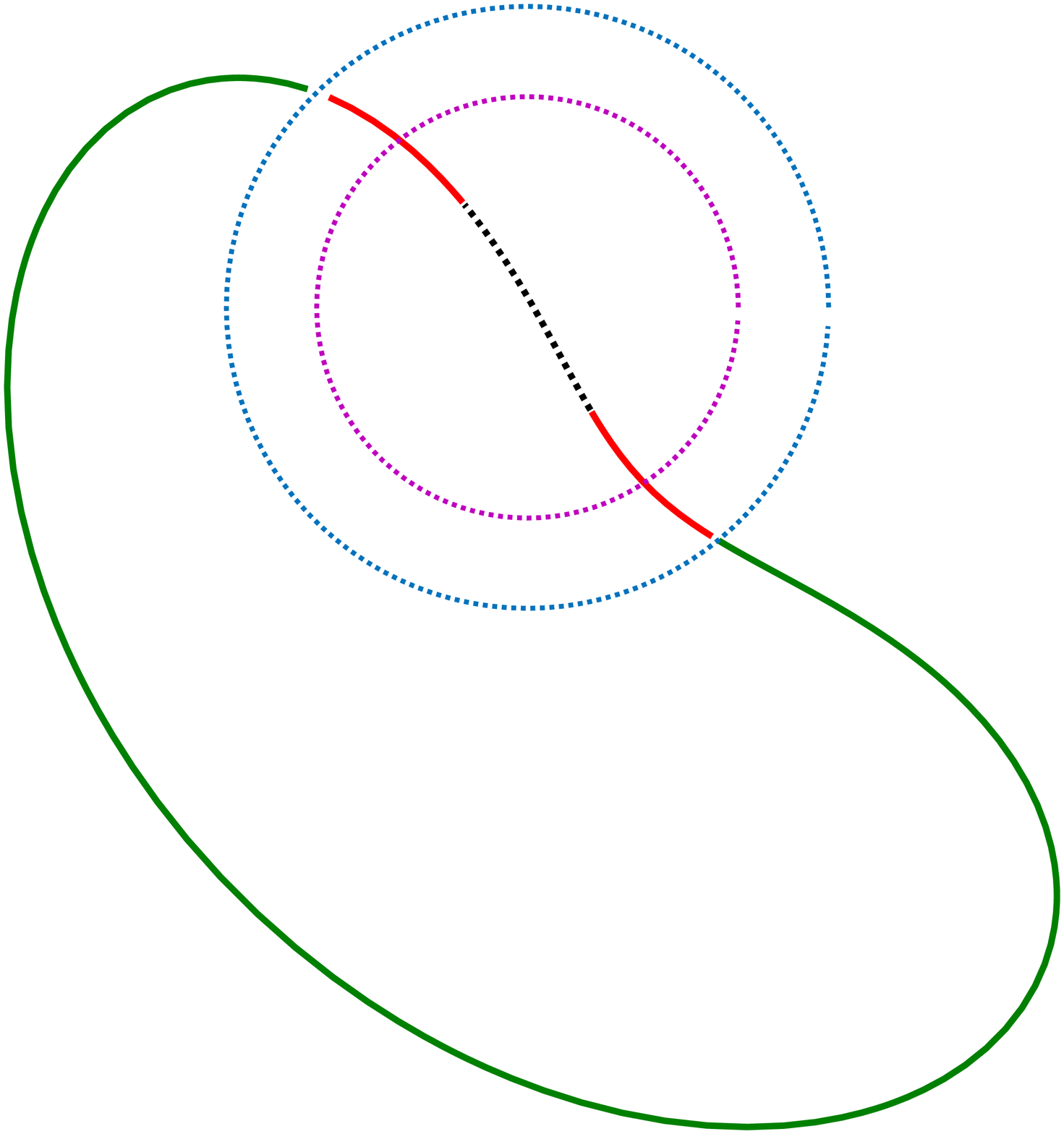}}
\put(110,95){$\Gamma_r$}
\end{picture}\\
(a) & &(b) 
\end{tabular}
\end{center}
\caption{ (a)  The proxy circle for $\Gamma_r$ shown in dash blue line divides $\Gamma_k$ into far (in green) and near (in red) with respect to $\Gamma_r$  
 (b) The interaction between the far-field part of $\Gamma_k$ and points on $\Gamma_r$
 can be captured by the interaction between points on $\Gamma_k^{\rm far}$ and a smaller proxy circle for  $\Gamma_r$ shown in dash purple. }
\label{fig:farNearSeparation}
\end{figure}

 The first step in creating the low rank factorization of $\mtx{Q}$ is constructing 
the low rank factorization of the non-zero blocks; i.e., the three factorizations in \agcmt{equation} (\ref{equ:block_lr}).
The construction of the low rank factorization of $\mtx{A}_{kp}$ starts with defining a circle 
$P^{\rm div}$ for $\Gamma_r$ which divides $\Gamma_k$ into two
 parts: the far-field and near-field with respect to $\Gamma_r$. Figure \ref{fig:farNearSeparation}(a) illustrates
 this separation. Let the superscript notation denote ``far'' and ``near'' parts of $\Gamma_k$. 
 \yzcmt{
 The separation corresponds to classifying the rows $\mtx{A}_{kp}$
 into two groups\agcmt{; the near- and far-field interactions}.   We first construct low-rank approximations to the far-field and near-field interaction separately and then merge them together for a final low-rank approximation $\mtx{L}_{kp}\mtx{R}_{kp}\approx \mtx{A}_{kp}$.}

 For far-field interaction, 
 the potential due to points in $\mathcal{I}_p$ evaluated at points on $\Gamma_k^{\rm far}$ can be approximated by 
 a linear combination of basis functions defined on any proxy surface that shields $\Gamma_r$ away from $\Gamma_k^{\rm far}$.
 Let $P^{\rm bas}$ denote the shielding proxy circle for $\Gamma_r$.  \agcmt{Here $P^{\rm bas}$ is 
 chosen to have a smaller radius as  $P^{\rm div}$ but the same center.  Figure \ref{fig:farNearSeparation}(b) illustrates
 an example of these circles.}  
 A low rank approximation for $\mtx{A}_{kp}^{\rm far}$ 
 can be constructed via an interpolative decomposition (ID) approximation \agcmt{(as defined below) for the matrix
 $\mtx{A}^{\rm far}_{k, \rm bas}$ which captures} the interaction between points on $\Gamma_k^{\rm far}$ and $P^{\rm bas}$. This is similar to the far-field compression idea in \cite{GILLMAN2012_china, MARPLE2016_periodicstokes, GILLMAN2013_quasiperiodic}.
 \yzcmt{The collection of skeleton row indices $J_k^{\rm far}$ from the ID for $\mtx{A}^{\rm far}_{k, \rm bas}$ correspond to discretiaztion points (or degrees of freedom) on $\Gamma_k$. Let $\mtx{P}_{k}^{\rm far}$ \agcmt{denote} the interpolation matrix, then a low-rank approximation to $\mtx{A}_{kp}^{\rm far}$  can be defined as $\mtx{A}_{kp}^{\rm far}\approx \mtx{L}_{kp}^{\rm far} \mtx{R}_{kp}^{\rm far}$ \agcmt{where} $ \mtx{L}_{kp}^{\rm far}= \mtx{P}_{k}^{\rm far}$ and $\mtx{R}_{k}^{\rm far}=\mtx{A}_{kp}^{\rm far}(J_k^{\rm far},:)$. 
 Here $\mtx{A}_{kp}^{\rm far}(J_k^{\rm far},:)$ denotes the submatrix of 
 $\mtx{A}_{kp}^{\rm far}$ with rows specified by $J_k^{\rm far}$.
 }

\begin{definition}
Given a tolerance $\epsilon$ and a $m\times n$ matrix $\mtx{W}$ (assuming $m<n$),
if there exists a positive integer $k_\epsilon\leq m$, and $m\times k_\epsilon$ 
matrix $\mtx{P}$ and vector $J$ such that
$$
\|\mtx{W} - \mtx{P}\mtx{W}(J(1:k_\epsilon),:)\|\leq \epsilon \|\mtx{W}\|,
$$
we call $\mtx{P}\mtx{W}(J(1:k_\epsilon),:)$ an interpolative decomposition (ID) approximation for $\mtx{W}$ with respect to the tolerance $\epsilon$.
Here $J$ is a vector of integers $j_i$ such that $1\leq j_i\leq m$ gives an ordering of the rows in $\mtx{W}$, and 
 $\mtx{W}(J(1:k_\epsilon),:)$ is a submatrix of $\mtx{W}$ with rows specified by the first $k_\epsilon$ entries of $J$.
$\mtx{P}$ is a $m\times k_\epsilon$ matrix that contains a $k_\epsilon \times k_\epsilon$ identity matrix.
Namely, $\mtx{P}(J(1:k_\epsilon),:) = \mtx{I}_{k_\epsilon}$.  The rows specified by $J(1:k_\epsilon)$ is called the \textit{skeleton row index}, and the matrix $\mtx{P}$ is referred to as the \textit{interpolation matrix}.

\end{definition}

\yzcmt{Due to the large number of discretization points on $\Gamma_k^{\rm far}$, it is often
too expensive to build the ID  for $\mtx{A}^{\rm far}_{k, \rm bas}$ directly.  %
Instead, we organize the discretization points on $\Gamma_k^{\rm far}$ into special structure such as the dyadic partition (See section 3 of \cite{Zhang_multilayer}) or binary tree (Such as the binary tree used in the HBS forward compression). 
The goal of using the special structure is to keep the cost of building the low-rank approximation linear with respect to the number of points on $\Gamma_k^{\rm far}$.
Then an ID for $\mtx{A}^{\rm far}_{k, \rm bas}$ is constructed by first building IDs for interaction between points in each individual partition subset or tree node and points on $P^{\rm bas}$, which corresponds to row subblocks of $\mtx{A}^{\rm far}_{k, \rm bas}$.  The individual IDs are then
combined into one final ID for $\mtx{A}^{\rm far}_{k, \rm bas}$. 
}
\yzcmt{\begin{remark}
Since the removed points $\mathcal{I}_c$ and added points $\mathcal{I}_p$ discretize the
same boundary curve segment $\Gamma_r$,  the far-field part of the low-rank approximation for $\mtx{A}_{kp}$ and $\mtx{A}_{kc}$ can be built from
the same ID approximation for $\mtx{A}^{\rm far}_{k, \rm bas}$.
The \agcmt{construction of the} approximations do not require explicit evaluation of the \agcmt{matrices}
$\mtx{A}_{kp}$ and $\mtx{A}_{kc}$.  \agcmt{Only the submatrices corresponding to the skeleton rows need to 
be evaluated for making the $\mtx{R}$ matrices.}
\end{remark}}

\begin{figure}[ht]
\begin{center}
\includegraphics[height=55mm,trim={3cm 0 3cm 0}, clip]{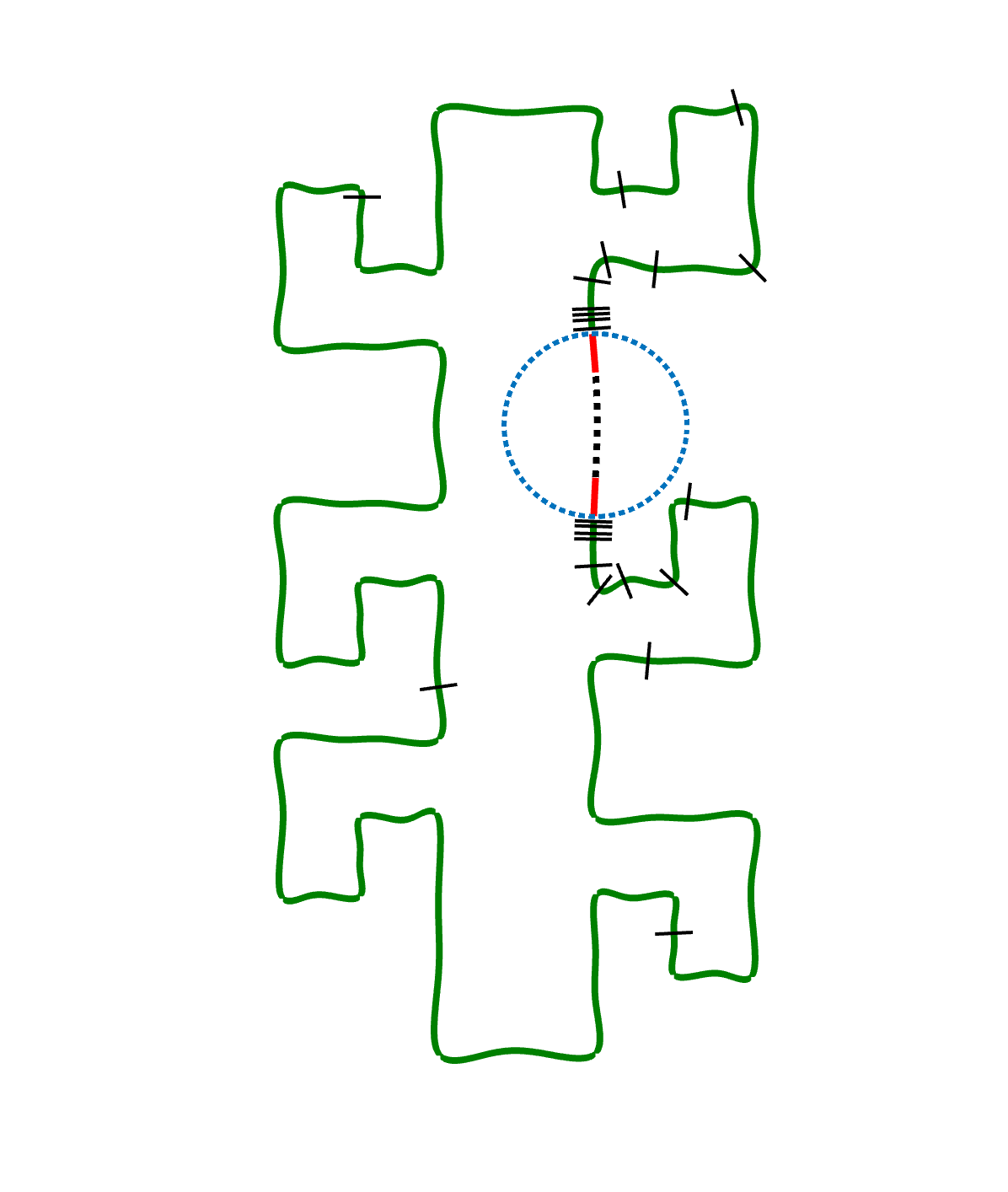}
\end{center}
\caption{A dyadic partition of $\Gamma_k^{\rm far}$ based on distance to $\Gamma_r$ for a refined channel discretization. The subintervals corresponding to a 10-level dyadic partition are plotted. }
\label{fig:ryan_pipe_dyadic}
\end{figure}

\begin{figure}[ht]
\begin{center}
\begin{tabular}{c c c c}
\includegraphics[height=45mm,trim={10cm 0 10cm 0}, clip]{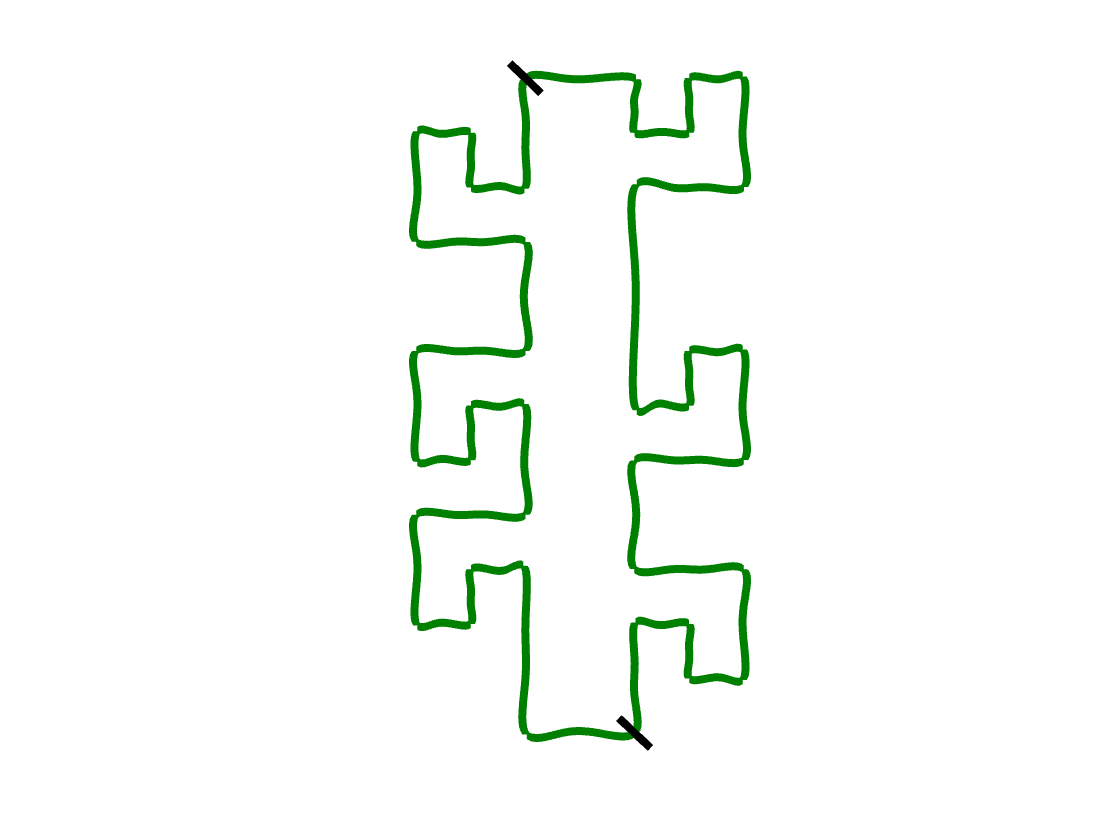}
 &
\includegraphics[height=45mm,trim={10cm 0 10cm 0}, clip]{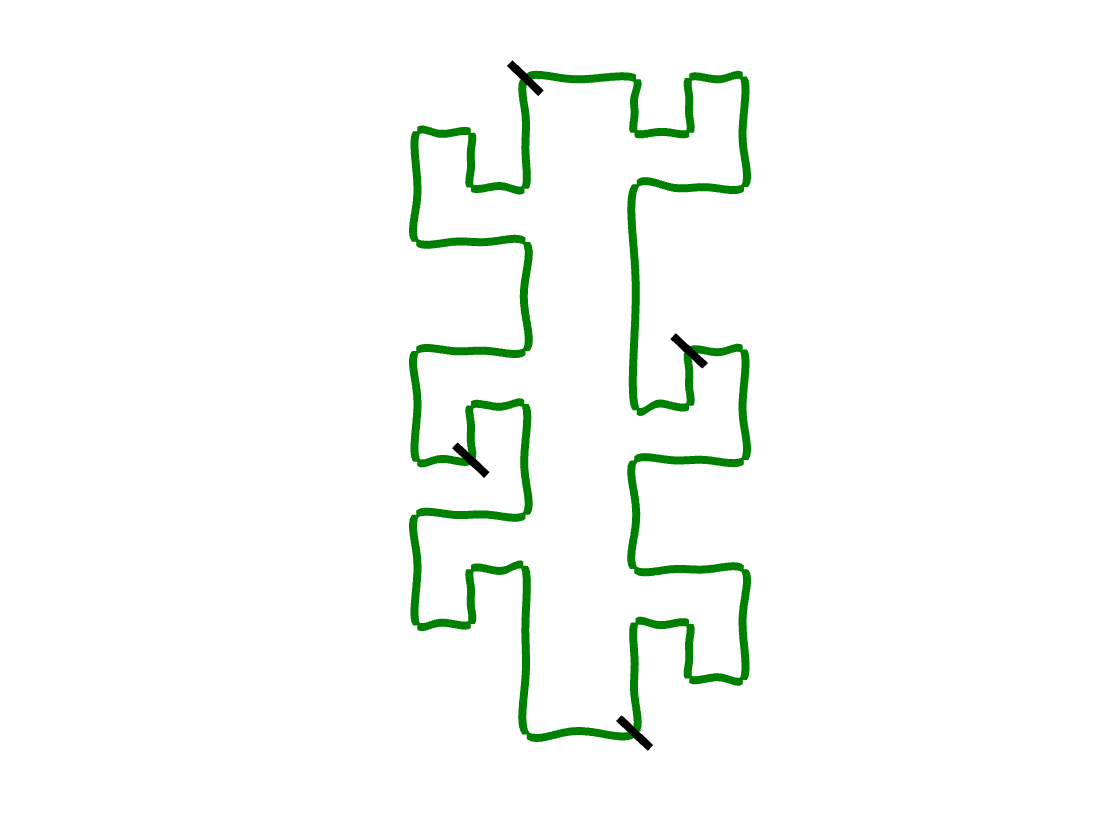}
 &
\includegraphics[height=45mm,trim={10cm 0 10cm 0}, clip]{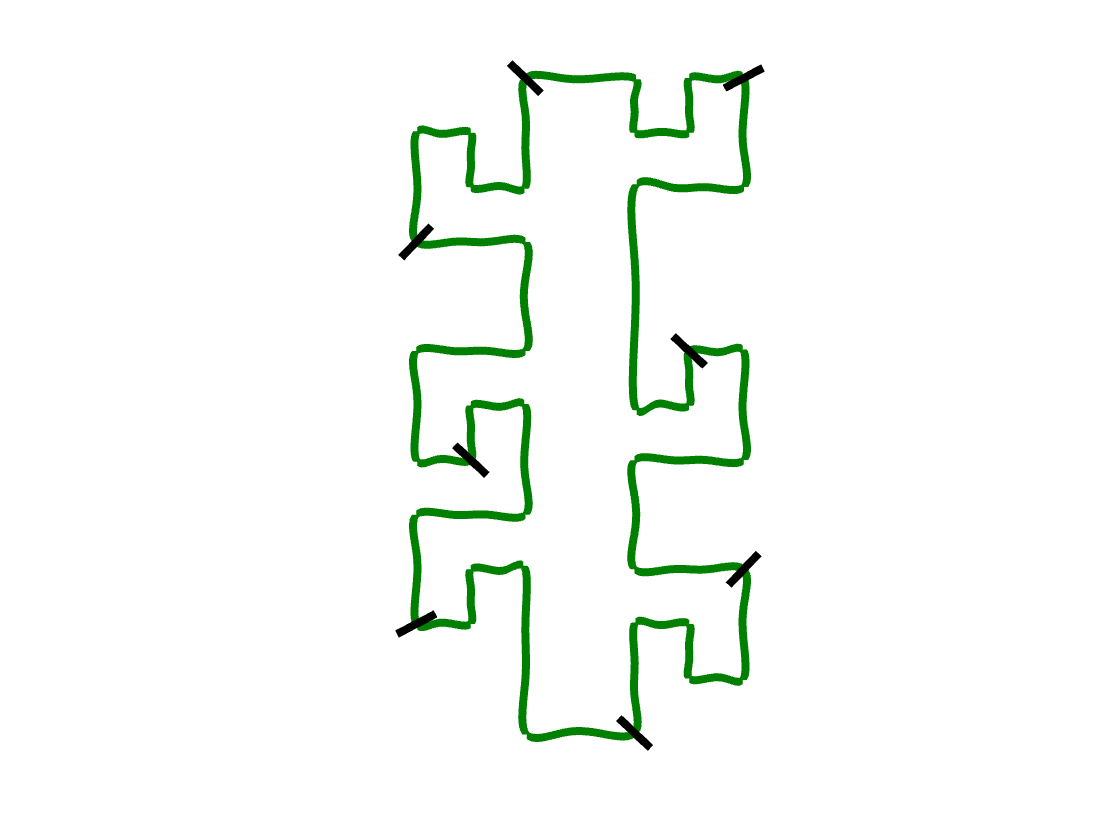}
&
\includegraphics[height=45mm,trim={10cm 0 10cm 0}, clip]{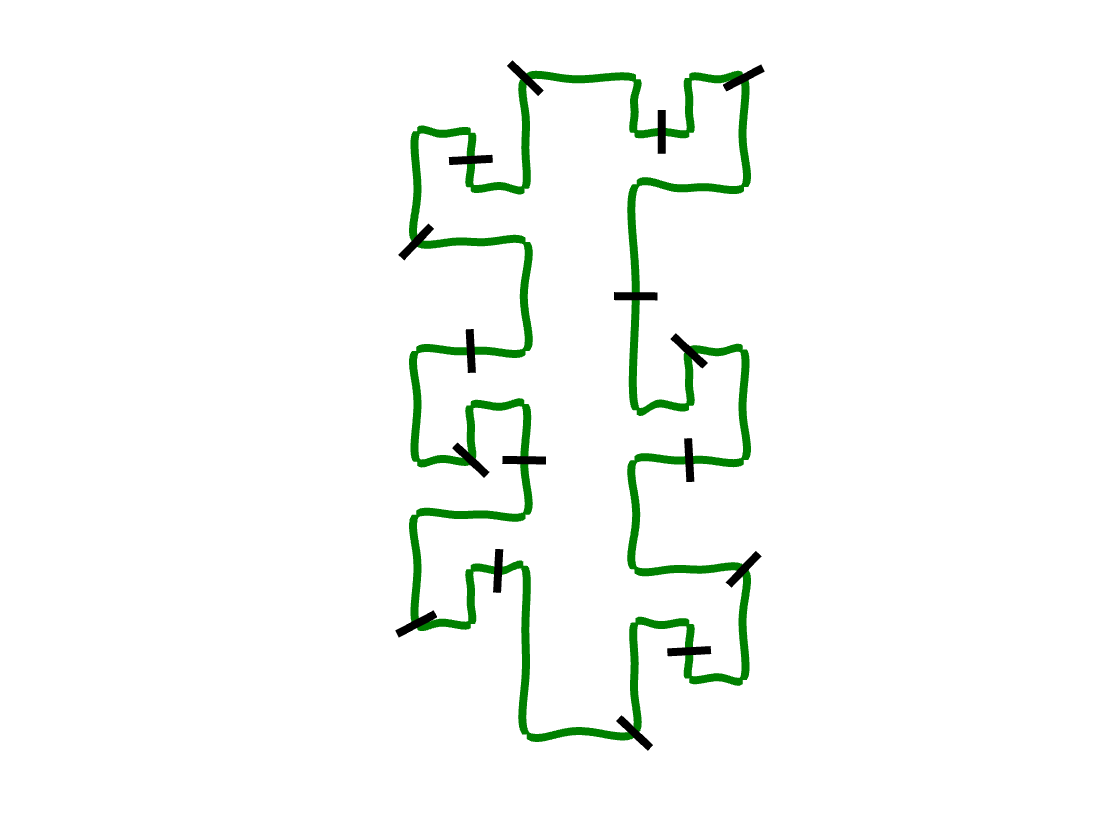}
 \\
Level 1 &Level 2&Level 3 &Level 4 \\
\end{tabular}
\end{center}
\caption{ Top four levels of a binary tree partition for the channel geometry. The subintervals (or boxes) corresponding to each level are plotted.  }
\label{fig:ryan_pipe_binary}
\end{figure}

\yzcmt{The choice of structure for creating the low rank factorization which will result in the
most efficient factorization technique depends
on how localized and the position of \agcmt{the portion of the boundary to be refined} $\Gamma_r$ relative to $\Gamma_k$.}
For example, the channel example given in section \ref{sec:scaling}
considers two kinds of local changes to the channel geometry in Figure \ref{fig:ryan_pipe}(a): 
\agcmt{a very localized refinement of the discretization illustrated in
Figure \ref{fig:ryan_pipe}(b); and a geometric perturbation consisting
of the addition of three interior circular holes as illustrated in Figure \ref{fig:ryan_pipe}(c).}
  For the problem in Figure \ref{fig:ryan_pipe}(b), 
the far-field and near-field separation is straightforward and
a dyadic partition of the far-field points on $\Gamma_k$ based on distance to $\Gamma_r$ is convenient and efficient.
However, for the problem in Figure \ref{fig:ryan_pipe}(c),
since the three holes do not cluster, 
a circle enclosing all holes would contain a large section of the channel boundary if not all of it, leading to lots of points on $\Gamma_k$ being clustered as ``near-field'' points although they are quite far away from any of the holes.
An efficient way to handle this problem is to introduce three circles each
enclosing an individual hole and define $P^{\rm bas}$ to be the union of the 
three circles.
And a binary tree, which does not \agcmt{have to } depend on distance to $\Gamma_r$, is 
a more appropriate choice.
Figure \ref{fig:ryan_pipe_dyadic} plots an example dyadic partition for the refined channel problem in Figure \ref{fig:ryan_pipe}(b), and  Figure \ref{fig:ryan_pipe_binary} plots the first three levels of an
example binary tree structure for the addition of holes problem in Figure \ref{fig:ryan_pipe}(c).

If there are not a large number of points that are near, 
which is often the case,
the near field interaction matrix $\mtx{A}_{kp}^{\rm near}$ can be compressed directly.
\yzcmt{
Otherwise, a dyadic partition of discretization points on $\Gamma_k^{\rm near}$ based on their distance to $\Gamma_r$ can be adopted.  The ID for $\mtx{A}_{kp}^{\rm near}$ 
can then be constructed in a hierarchical way utilizing the idea of tree-node wise proxy circles (See section 3 of \cite{Zhang_multilayer}).}

Once both far and near part of $\mtx{A}_{kp}$ are compressed, the low rank factors can be concatenated 
to form a low rank approximation for  $\mtx{A}_{kp}$.
One may want to apply ID again to the concatenated factors to further reduce the rank numbers.
%
%
 
The near-field part of  $\mtx{A}_{kc}$ can be constructed in similar way as that of $\mtx{A}_{kp}$.
For $\mtx{A}_{pk}$, we consider again a far-field and near-field separation of the points on $\Gamma_k$ based on distance to $\Gamma_r$, which corresponds to classifying the columns of the matrix into two groups.
The far-field interaction $\mtx{A}_{pk}^{\rm far}$ can be obtained by an ID of $\mtx{A}_{p,\rm div}$ the interaction between the added points discretizing $\Gamma_r$ and sample points on the separation circle $P^{\rm div}$.
 If the number of points added is large,
we can relieve the computational burden by using a dyadic partition or binary tree as for building the ID for $\mtx{A}_{k,\rm bas}^{\rm far}$.
The construction for $\mtx{A}_{pk}^{\rm near}$ is similar to the near-field part of the approximation for the near-field of  $\mtx{A}_{kp}$ and $\mtx{A}_{kc}$.
%
\yzcmt{
\begin{remark}
When approximating the three blocks in $\mtx{Q}$, we always use ID to compress the rows of the matrices. We also uniformly define the $\mtx{L}$ factor of the low-rank approximation to be the interpolation matrix (or product of multiple interpolation matrices if special tree structure is used) and the $\mtx{R}$ factor to be the submatrix of the discretized BIE specified by the skeleton row indices given by the IDs. This uniform format for all three blocks is intentional as it improves the conditioning of applying the Woodbury formula.
More details will be given in Section \ref{sec:woodbury_stability}. 
Note the blockwise compression technique given in \cite{ZHANG2018} manages to compress all far-field part of the three blocks $\mtx{A}_{kp}$, $\mtx{A}_{kc}$, and $\mtx{A}_{pk}$ using one binary tree by doing row-wise ID for $\mtx{A}_{kp}$ and $\mtx{A}_{kc}$ but column-wise ID for 
$\mtx{A}_{pk}$. Namely, the far-field for all three blocks are approximated by the same set of skeleton points on $\Gamma_k$. 
For Laplace problems, the technique in \cite{ZHANG2018} is expected to be more efficient than the one presented here especially for the case where $\Gamma_k^{\rm far}$ contains lots of points.
But for Stokes problems, the mixed usage of row- and column-wise ID leads to 
conditioning issues and should be avoided. 
\end{remark}}
\yzcmt{With the special structure and partitioning, the cost of constructing the low-rank factorization for $\mtx{A}_{kp}$ is $O\left((N_k+N_p)k_{kp} \right)$.
Similarly, the cost for factorizing $\mtx{A}_{kc}$ is $O\left((N_k+N_c)k_{kc}\right)$. 
And the cost of factoring $\mtx{A}_{pk}$ is $O\left((N_k+N_p)k_{pk}\right)$.
}

\subsection{Stable application of the Woodbury formula}
\label{sec:woodbury_stability}
Woodbury formulas such as (\ref{eq:wood}) are well-known in the linear algebra literature \cite{GOLUB1996} and have been the cornerstone of recently developed fast direct solvers
for applications including  periodic Stokes flow \cite{MARPLE2016_periodicstokes} and 
quasi-periodic scattering problems \cite{GILLMAN2013_quasiperiodic,Zhang_multilayer}. 
While the Woodbury formulas have been used in these applications, it was done so without
any concern for the stability of the approach.  
\yzcmt{
This section will review the 
stability analysis of the Woodbury formula given in \cite{YIP1986},
investigate its use in the case of 
Stokes problems,
and \agcmt{presents} the two-step construction of the low-rank factorization of $\mtx{Q}$
started in the previous section.
}

The main concern in the stability of the Woodbury formula \agcmt{lies
in the stable inversion of the matrix $\mtx{W} = \mtx{I}+\mtx{R}\tilde{\mtx{A}}^{-1}\mtx{L}$.  We will refer to the matrix $\mtx{W}$ as the \textit{Woodbury operator}.}
\cite{YIP1986} states that in order to stably solve a linear system via the Sherman-Morrison-Woodbury formula,  the following two conditions must be satisfied by the linear system:
\begin{itemize}
 \item[(i)]
  All the relevant matrix-matrix and matrix-vector multiplications in (\ref{eq:wood}) involving $\tilde{\mtx{A}}^{-1}$ are numerically stable.
  \item[(ii)] The Woodbury operator is well-conditioned.
  \end{itemize}
  For Stokes problems, the first condition is satisfied thanks to the choice of boundary integral formulation (in Section \ref{sec:BIE}) 
  and the use of a stable fast direct solver. 
  \agcmt{Since Stokes problems tend to have a large condition number, we choose to modify the second condition to: (ii) The Woodbury 
  operator is as ``well-conditioned" as the full linear system $\tilde{\mtx{A}}+\mtx{L}\mtx{R}$.}

\agcmt{The following lemma, which is a modified version of Lemma 1 in \cite{YIP1986},}  \yzcmt{provides an upper bound on the condition number of the Woodbury \agcmt{operator and formally defines what we mean by the Woodbury operator being as ``well-conditioned" as the fully linear system.}
\agcmt{The lemma is stated} in the context of discretized boundary integral operators and the ELS.}
Specifically it
provides conditions on the low rank approximation of the update matrix $\mtx{Q} \approx \mtx{LR}$
which must be satisfied (along with both the linear systems for the original and refined
discretization being well-conditioned) \agcmt{for the Woodbury operator to be ``well-conditioned."}

\begin{lemma}[Upper bound on the condition number of the Woodbury operator]
\label{thm:woodbury}
 Assume the operator $\mtx{A}_{oo}$, $\mtx{A}_{nn}$, $\mtx{A}_{cc}$ and $\mtx{A}_{pp}$ as defined in (\ref{equ:originalLinearSystem}), (\ref{equ:perturbedlLinearSystem})  and (\ref{eq:extendedlLinearSystem}) are all invertible. 
 \yzcmt{Then the operator  $\mtx{A}_{\rm ext}$ in the ELS $\mtx{A}_{\rm ext}\boldsymbol{\sigma}_{\rm ext}=\mtx{g}_{\rm ext}$  is invertible.}
  Let $\hat{\mtx{A}}_{\rm ext}= \tilde{\mtx{A}}+\mtx{LR}$ denote the approximation of $\mtx{A}_{\rm ext}$.
   If the $k$ columns in the low-rank factor $\mtx{L}$ and
   the $k$  rows in $\mtx{R}$ are linearly independent, then the 
   condition number of the Woodbury operator is bounded above as follows:

      \begin{equation}\label{equ:woodburyConditionThm}
 \begin{split}
 \kappa\left( \mtx{I} +\mtx{R}\tilde{\mtx{A}}^{-1}\mtx{L} \right)
&\leq \min\left\{ \hat{\kappa}(\mtx{L})^2 , \,    \hat{\kappa}(\mtx{R})^2   \right\} \kappa\left( \hat{\mtx{A}}_{\rm ext}\right)  \kappa\left(\tilde{\mtx{A}} \right),
\end{split}
 \end{equation}
 where
 $$  \hat{\kappa}(\mtx{L})= \left \|   \mtx{L}^\dag\right\|  \left \|   \mtx{L}\right\| 
 \mbox{ and } 
 \hat{\kappa}(\mtx{R})= \left \|   \mtx{R}^\dag\right\|  \left \|   \mtx{R}\right\|
 $$
 with
 $$
  \mtx{L}^\dag =\left( \mtx{L}^T \mtx{L}\right)^{-1}\mtx{L}^T
 \mbox{ and }
  \mtx{R}^\dag =\mtx{R}^T\left( \mtx{R} \mtx{R}^T\right)^{-1}
 $$
 defined as the pseudo-inverse for $\mtx{L}$ and $\mtx{R}$ in the standard sense.
 \label{thm:bound}
 \end{lemma}
The proof of the lemma can be found in \cite{YIP1986} and is also included in the Appendix.
 
\begin{figure}[ht]
\begin{center}
\includegraphics[height=45mm, trim={0 0.5cm 0 0.5cm}]{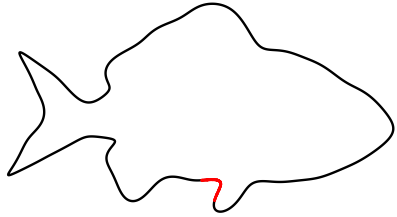}

\end{center}
\caption{Illustration of a fish geometry where the red portion of the boundary is 
refined. Table \ref{tab:conditiontest} reports on the conditioning of the Woodbury 
operator for different numbers of discretization points and refinements.
}
\label{fig:conditiontest}
\end{figure}

Thus, if we construct $\mtx{Q}\approx \mtx{LR}$ so that  $\hat{\mtx{A}}_{\rm ext} =\tilde{\mtx{A}}+\mtx{LR}$ is invertible, 
 $\mtx{L}$ and $\mtx{R}$ are full-rank, and additionally let  $ a^2 = \min\left\{ \hat{\kappa}(\mtx{L})^2 , \,    \hat{\kappa}(\mtx{R})^2   \right\}$, then $\kappa\left( \mtx{I} +\mtx{R}\tilde{\mtx{A}}^{-1}\mtx{L} \right)
\leq a^2\,\kappa\left( \hat{\mtx{A}}_{\rm ext}\right)  \kappa\left(\tilde{\mtx{A}} \right)$. 
When the original problem and new problem have similar condition numbers, i.e., $  \kappa\left( \hat{\mtx{A}}_{\rm ext}\right) \approx \kappa\left(\tilde{\mtx{A}} \right) \approx \kappa$, 
the lemma and the low-rank approximation construction above 
together give the bound $\kappa\left( \mtx{I} +\mtx{R}\tilde{\mtx{A}}^{-1}\mtx{L} \right)\leq a^2 \kappa^2$. 
\yzcmt{The upper bound given by the lemma can be improved by building $\mtx{L}$ and $\mtx{R}$
so that at least one of $\hat{\kappa}(\mtx{L})$ and $\hat{\kappa}(\mtx{R})$ stay small.
One way to do this for the update matrix $\mtx{Q}$ is to build a truncated SVD for $\mtx{Q}$ for some given tolerance and assign $\mtx{L}$ to be the semi-unitary matrix corresponding to the column space and $\mtx{R}$ to be the rest of the factors in the decomposition.
By doing this, matrix $(\mtx{L}^T\mtx{L})$ and $(\mtx{R}^T\mtx{R})$ stays away from being singular,
leading to minimal values of $\hat{\kappa}(\mtx{L})$ and $\hat{\kappa}(\mtx{R})$.
Since $a^2$ is defined to be the smaller one among
$\hat{\kappa}(\mtx{L})^2$ and $\hat{\kappa}(\mtx{R})^2$ \agcmt{in \eqref{equ:woodburyConditionThm}},
only one of  $\hat{\kappa}(\mtx{L})^2$ and $\hat{\kappa}(\mtx{R})^2$ being small
is sufficient.
For example, if we construct the truncated SVD for each of the non-zero blocks in $\mtx{Q}$,
$$\mtx{U}_{pk}\mtx{\Sigma}_{pk}\mtx{V}^T_{pk}\approx \mtx{A}_{pk},\;
\mtx{U}_{kc}\mtx{\Sigma}_{kc}\mtx{V}^T_{kc}\approx \mtx{A}_{kc},\mbox{ and }
\mtx{U}_{kp}\mtx{\Sigma}_{kp}\mtx{V}^T_{kp}\approx \mtx{A}_{kp}$$
then we would define the concatenated factors for the updated matrix $\mtx{L}_{\rm block}$ and $\mtx{R}_{\rm block}$ so that 
the three non-zero blocks are uniform in format: for example,
\begin{equation}\label{equ:uniform_concat}
\mathbf{L}_{\rm block}=
\begin{bmatrix}
\mtx{0} & -\mtx{U}_{kc} & \mtx{U}_{kp}\\
\mtx{0} &\mtx{0} & \mtx{0}\\
\mtx{U}_{pk} &\mtx{0}&\mtx{0}\\
\end{bmatrix}  
\mbox{ and }
\mathbf{R}_{\rm block}=
\begin{bmatrix}
\mtx{\Sigma}_{pk}\mtx{V}^T_{pk} & \mtx{0} &\mtx{0}\\
\mtx{0} & \mtx{\Sigma}_{kc}\mtx{V}^T_{kc} &\mtx{0}\\
\mtx{0}  & \mtx{0} & \mtx{\Sigma}_{kp}\mtx{V}^T_{kp}
\end{bmatrix}.
\end{equation}
And as a final step to push the factorization closer to optimal,
one would build a truncated SVD for $\mtx{L}_{\rm block}\approx \mtx{U}_{\rm block}\mtx{\Sigma}_{\rm block}\mtx{V}^T_{\rm block}$ and define the finalized low-rank approximation as
$\mtx{Q}\approx \mtx{L}_{\rm optimal}\mtx{R}_{\rm optimal}$  with $\mtx{L}_{\rm optimal}=\mtx{U}_{\rm block}$ and $\mtx{R}_{\rm optimal}=\mtx{\Sigma}_{\rm block}\mtx{V}^T_{\rm block}\mtx{R}_{\rm block}$.
As a demonstration, 
 Table \ref{tab:conditiontest} illustrates the conditioning of the Woodbury operator corresponding to the concatenated factorization $\mtx{Q}\approx \mtx{L}_{\rm block}\mtx{R}_{\rm block}$ and  the factorization after the extra SVD refactorization $\mtx{Q}\approx \mtx{L}_{\rm optimal}\mtx{R}_{\rm optimal}$ for  
 a sample problem defined on the fish geometry illustrated in Figure \ref{fig:conditiontest}. 
The tolerance for the SVD truncation is set to $10^{-10}$ and the condition numbers reported in the table are 
 calculated via Matlab's \texttt{cond()} function.
}

\begin{table}[ht]
\begin{center}
\begin{tabular}{|c||c|c||c|c||c|c||c|}
\hline
$N_k,N_c,N_p$ & $k_{\rm block}$ & $\kappa_{\rm block}$ & $k_{\rm optimal}$& $\kappa_{\rm optimal}$&$\kappa\left( \hat{\mtx{A}}_{\rm ext}\right)$ & $\kappa\left(\tilde{\mtx{A}} \right)$& Upper bound\\
\hline
752, 48, 384 &175 &1630.0 & 141 &  98.5 & 378.0 & 371.0 &4.1e+25\\
1520, 80, 640 &158 &407.3 & 124  &  77.3  & 376.7 & 371.0 &1.2e+25\\
3072, 128, 1024 &143 &523.5 &  113 &  77.0 & 375.0 & 371.0 &1.6e+25\\
\hline
\end{tabular}
\end{center}
\caption{The observed rank (same as the size of the Woodbury system) and condition number for an interior BIE on the fish geometry illustrated in Figure \ref{tab:conditiontest}. 
$k_{\rm block}$ and $\kappa_{\rm block}$ are the size and condition number for the Woodbury operator defined for the concatenated factorization $\mtx{Q}\approx \mtx{L}_{\rm block}\mtx{R}_{\rm block}$.
$k_{\rm optimal}$ and $\kappa_{\rm optimal}$ are the size and condition number for the Woodbury operator defined for the factorization with an extra SVD applied to $\mtx{L}_{\rm block}$, i.e., $\mtx{Q}\approx \mtx{L}_{\rm optimal}\mtx{R}_{\rm optimal}$.
All the factorizations are constructed by truncating SVD to the desired accuracy of $10^{-10}$. The condition numbers are all calculated by Matlab's \texttt{cond()} function.
  The condition number of the block-diagonal matrix $\tilde{\mtx{A}}$ and the 
 ELS ${\mtx{A}}_{\rm ext}$ are also reported.  Finally the upper bound given by 
 Lemma \ref{thm:bound} corresponding to $\mtx{Q}\approx \mtx{L}_{\rm optimal}\mtx{R}_{\rm optimal}$ is also provided.}
\label{tab:conditiontest}
\end{table}

\yzcmt{
Truncated SVDs are expensive to construct. \agcmt{Therefore} this
 optimal approach is not computationally viable except for problems small in size.
Instead we propose an alternative two-step approach which \agcmt{addresses} the above conditioning considerations but is less expensive and thus suitable for large size problems.
The first step is to construct the 
 low-rank approximations for block $\mtx{A}_{pk}$, $\mtx{A}_{kc}$ and $\mtx{A}_{kp}$
following the method given in section \ref{sec:low_rank_approx}. 
While this constructs a valid low rank factorization of $\mtx{Q}$, the approximation is 
often quite far away from being optimal and results in an unnecessarily large condition number of  the Woodbury operator. 
\agcmt{In fact, there are many cases where the resulting Woodbury operator 
is ill-conditioned even though the original system is well-conditioned.}
To remedy this \agcmt{artifical poor conditioning, we propose 
the refactorization of $\mtx{L}_1$ via the} random sampling based ID decomposition.  \agcmt{This compresses the rows of $\mtx{L}_1$ and results
in a significantly closer to optimal rank factorization. It is 
important to maintain a uniform format in the refactorization technique by 
always} applying \agcmt{the} ID to compress the rows of
matrices when building the blockwise \agcmt{factorization} and assigning the interpolation matrix factor from the ID approximation to be blocks in $\mtx{L}_1$.
}

\yzcmt{
\begin{remark}
For the problems considered in this manuscript, \agcmt{the upper bound in Lemma \ref{thm:woodbury}} is  overly pessimistic. \agcmt{In fact,} the observed condition number is much smaller than $a^2\kappa^2$. 
\agcmt{In practice, when }the low-rank approximation \agcmt{of} $\mtx{Q}$ is constructed with care \agcmt{via the SVD technique or via the  two-step compression based on IDs}, the observed condition number of the Woodbury system \agcmt{is} comparable to the condition number of the original linear system. 
 For example,
 Table \ref{tab:conditiontest}  also reports the upper bound on the condition number  
 given by Lemma \ref{thm:bound} for the Woodbury system with the final factorization.
 While rank and condition number are improved by the extra SVD
 recompression, both condition numbers are well below the upper bound provided by the lemma.
 \end{remark}
}

\yzcmt{The cost for the extra ID refactorization is $O((N_p+N_k+N_c)k_1k)$. Thus, the
total cost for constructing the final low-rank approximation for the update matrix scales
linearly with respect to $N_k+N_p+N_c$.
}

\section{A preconditioner for BVPs on locally refined discretization}
\label{sec:preconditioner}
The ELS presented in Section \ref{sec:localperturb_ext} is very useful for 
problems where there is local refinement \agcmt{of the} discretization.  While the fast direct solver for the ELS is efficient, it can suffer from a loss in accuracy when the problem  \agcmt{has a high condition number. This is frequent occurrence for Stokes problems especially in complex geometries.}  An alternative \agcmt{to fast direct solvers} is to use an iterative solver 
coupled with a fast matrix vector multiplier such as the FMM in these instances.  \agcmt{The large condition number often means that a large number of iterations are required for the iterative solver to converge.}  This section 
presents an alternative \agcmt{solution technique which is essentially the union of a fast direct solver with an iterative solver.}  Roughly speaking, the 
technique is to use the direct solver presented in Algorithm \ref{alg:fdsolver} 
as preconditioner for the ELS that is solved via an iterative solver  
coupled with a fast matrix vector multiplier.  

Section \ref{sec:general_hbs_preconditioner} details how the accuracy in which the direct 
solver is constructed impacts its ability to be a preconditioner.  Then Section \ref{sec:local_perturbation_preconditioner}
details the preconditioner developed for the ELS (\ref{eq:extendedlLinearSystem}).  


\subsection{HBS inverse approximation as preconditioner}
\label{sec:general_hbs_preconditioner}
It is becoming more common to use low accuracy fast direct solvers as preconditioners for linear 
systems that arise from discretizations of integral equations and differential equations \cite{Darve_IFMM_pre,Hmat_pre,2010_HSS_pre,2005_beb,2003_beb}.  This section explores effectiveness of fast direct solvers as preconditoners
for \agcmt{the discretized} integral equation \agcmt{associated with} an interior Stokes problem.  

Consider the linear system $\mtx{A}\bm{\sigma}=\vct{g}$ which results from the discretization of \agcmt{equation} (\ref{equ:interior_BIE}).   
 Let $\epsilon$ denote the tolerance for which the fast direct solver  was constructed
 and $\mtx{A}_{\epsilon}^{inv}$ denote the corresponding approximate inverse
 of $\mtx{A}$.  Then the left-preconditioned problem is defined as
\begin{equation}\label{eq:left_precond_general}
\left(\mtx{A}_{\epsilon}^{inv}\mtx{A}\right)\bm{\sigma}=\left(\mtx{A}_{\epsilon}^{inv}\vct{g}\right).
\end{equation}

To investigate the performance of the fast direct solver as a preconditioner with different tolerances $\epsilon$,
we consider the fish geometry in Figure \ref{fig:conditiontest} 
with no local refinements. 
In particular, we place two hundred 16-point Gaussian panels uniform in parameterization space on the boundary. 
 The linear system
(\ref{eq:left_precond_general}) is solved via GMRES\cite{SaadSchultz1986}.  \agcmt{The application of $\mtx{A}$
and $\mtx{A}_\epsilon^{inv}$ is done via the HBS technique from \cite{GILLMAN2012_china}. The performance of the solver will be the 
same for any fast direct solver.}  The tolerance for the compression of the matrix vector operator is fixed at
$10^{-10}$.  The time for constructing the HBS \agcmt{representation of the} matrix is $6.81$ seconds on a single core 1.6GHz 8GB RAM desktop.
Table \ref{tab:precondition_on_fish} reports the performance of the preconditioned 
solution technique.  For all experiments, the tolerance of the iterative solver is set to $10^{-11}$ \agcmt{and the average relative error  in the solution compared against the exact solution at sampled interior locations is roughly $7\times 10^{-10}$}.   
Recall \agcmt{from Table \ref{tab:conditiontest}} that the linear system is well conditioned.
Thus even without a preconditioner, only 55 iterations are needed  to achieve the desired tolerance.  
The results indicate that \agcmt{low accuracy approximations ($\epsilon>10^{-3}$) do not improve the performance of the iterative 
solver enough to justify constructing the preconditioner.}  \agcmt{For 
$\epsilon <10^{-3}$,} the minimum number of repeated solves needed to justify the use 
of the preconditioner grows as $\epsilon$ decreases. 
\agcmt{This experiment illustrates that the use of a low accuracy 
fast direct solver as preconditioner is not fruitful in improving the 
convergence rate of iterative solvers.}
For \agcmt{problems where the condition number of the discretized linear system is large, a preconditioner} may be required for the iterative solver to converge within reasonable number of iterations \agcmt{with the available} computing resources.

\begin{table}[h]
\begin{center}
\begin{tabular}{|c|c|c|c|c|}
\hline
$\epsilon$& $n_{\rm iter}$ & $T_{\rm pre}$ & $T_{\rm sol}$ & MinSol \\
\hline
No preconditioner & 55 & NA & 5.2e-1 &NA\\
\hline
 1e-10 & 2 & 7.66 &6.6e-2&17\\
\hline
1e-8 & 2&4.92 &1.0e-1&12\\
\hline
 1e-6 & 4& 2.89&1.3e-1&8\\
\hline
1e-5 & 6&2.23 &1.5e-1&7\\
\hline
 1e-4 & 11&1.74 &2.2e-1&6\\
\hline
 1e-3 & 36& 1.12&5.1e-1&--\\
\hline
 1e-2 &52 & 1.11&7.9e-1&--\\
\hline
 1e-1 &53 & 0.97&8.1e-1&--\\
\hline
\end{tabular}
\caption{Number of iteration $n_{\rm iter}$, time in seconds to build the preconditioner $T_{\rm pre}$, 
time in seconds for GMRES to converge $T_{\rm sol}$ and the minimum number of solves MinSol needed to 
justify the use of the preconditioner when using an HBS inverse approximation with accuracy $\epsilon$  as 
a preconditioner for the interior BIE on the fish geometry in Figure \ref{fig:conditiontest}. 
The boundary geometry is discretized with two hundred 16-point Gaussian panels uniformly distributed in parameterization space.
With this discretization, the average relative solution error at sample locations on the interior is roughly $7\times 10^{-10}$.
}
\label{tab:precondition_on_fish}
\end{center}
\end{table}

\subsection{Preconditioned iterative solver for the locally refined problem}
\label{sec:local_perturbation_preconditioner}
Just like the discreitzed BIE for a Stokes boundary value problem on a given geometry, the ELS (\ref{eq:extendedlLinearSystem}) can also suffer from conditioning issues. This section presents a preconditioner
based on the solver from Section \ref{sec:fdsolver} and a fast matrix vector multiplier that can be utilized
to accelerate an iterative solver.  It is expected that the number of iterations needed to converge will be less than if there was 
no preconditioner at all.  Additionally, there is no loss of digits associated with inverting \agcmt{poorly} conditioned matrices.

The idea behind the preconditioner is simple.  Let $\mtx{A}^{inv}_{oo}$ and $\mtx{A}^{inv}_{pp}$ denote the 
approximate (or exact if the matrices are small enough) inverses of $\mtx{A}_{oo}$ and $\mtx{A}_{pp}$,
respectively.  Then 

$$\tilde{\mtx{A}}^{inv}
= \begin{bmatrix}
\mtx{A}_{oo}^{inv} & \mtx{0} \\
 \mtx{0} & \mtx{A}_{pp}^{inv}\\
\end{bmatrix}
\approx
 \tilde{\mtx{A}}^{-1}
$$
and
$$
\mtx{A}_{\rm ext}^{inv}=
\tilde{\mtx{A}}^{inv}
-\tilde{\mtx{A}}^{inv}\mtx{L}
  \left( \mtx{I} + \mtx{R}\tilde{\mtx{A}}^{inv1}\mtx{L}\right )^{-1}
  \mtx{R}\tilde{\mtx{A}}^{inv}
  \approx \mtx{A}_{\rm ext}^{-1}.
$$
The Woodbury formula can be applied efficiently to any vector via the technique presented in Algorithm \ref{alg:fdsolver}.  

Instead of solving the true ELS, we propose solving the approximation of the linear system (\ref{eq:extendedlLinearSystem}) where $\mtx{A}_{\rm ext}$ is approximated by a block diagonal plus low rank form; 
i.e., as $\mtx{A}_{\rm ext} \approx \left( \tilde{\mtx{A}} +\mtx{LR}\right)$.
The matrix $\tilde{\mtx{A}}$ can be applied to a vector $\vct{b}$ block-wise
 $$\tilde{\mtx{A}}\vct{u}=\begin{bmatrix}
\mtx{A}_{oo} & \mtx{0} \\
 \mtx{0} & \mtx{A}_{pp}\\
\end{bmatrix}\begin{bmatrix}
\vct{b}_o\\
\vct{b}_p
\end{bmatrix}
=\begin{bmatrix}
\mtx{A}_{oo}\vct{b}_o\\
\mtx{A}_{pp}\vct{b}_p
\end{bmatrix}.
$$
The evaluation of $\mtx{A}_{oo}\vct{b}_o$ can be accelerated via fast matrix-vector multiplication algorithms, such as \agcmt{the} FMM or \agcmt{the approximate forward operator created in the process of building 
a fast direct solver, and is constructed for the original discretization.}
Similar to the fast direct solver for the ELS presented in Algorithm \ref{alg:fdsolver},
if $N_p$ is small, \agcmt{the matrix} $\mtx{A}_{pp}$ can be constructed and applied via dense linear algebra.
Otherwise, a separate fast matrix-vector multiplication can be constructed for $\mtx{A}_{pp}$.
Since $\mtx{L}$ and $\mtx{R}$ are block sparse  and low-rank, 
they can be applied to any vector densely with little cost.

In this paper, we assume a forward HBS representation, the HBS inverse, and matrix-vector multiplication for applying $\mtx{A}_{oo}$ and $\mtx{A}_{oo}^{inv}$ are available. Then the ELS  for the problem defined on the refined geometry 
\begin{equation}\label{eq:ELS_approximate}
\left( \tilde{\mtx{A}} +\mtx{Q}\right)\bm{\tau}_{\rm ext}\approx\left( \tilde{\mtx{A}} +\mtx{LR}\right)\bm{\tau}_{\rm ext}= \vct{g}_{\rm ext}    
\end{equation} \agcmt{only requires building the low rank factorization of the blocks in $\mtx{Q}$ and the operators associated with the $\mtx{A}_{pp}$ block and} \yzcmt{can be solved by an iterative solver such as GMRES. $\mtx{A}_{\rm ext}^{inv}$ can be constructed with the extra cost of carrying out the Woodbury formula and applied as a preconditioner to  (\ref{eq:ELS_approximate}).}
\agcmt{For a well-conditioned problem, where many different choices of local refinements
and/or right-sides are considered,}
the total cost may be greatly reduced by using the fast direct solver in Section \ref{sec:fdsolver} as a preconditioner.   Table \ref{tab:precondition_on_refinedfish} reports the performance of the 
 preconditioner when it is applied to the boundary value problem on the fish geometry in Figure \ref{fig:conditiontest} 
 where
 the red region of the boundary is refined.  
 \yzcmt{The original discretization has two hundred 16-point Gaussian panels uniformly distributed in parameterization space;
 8 panels discretize the red region and are replaced by 64 panels for the refinement.}
\agcmt{The number of discretization points
kept was $N_k = 3072$, the number of discretization points cut was $N_c = 128$ and 
the number of discretization points added was $N_p=1024$.}
The tolerance for HBS compression and low-rank approximations \agcmt{were} set to $10^{-10}$, and the tolerance for GMRES \agcmt{was} set to $10^{-11}$. The average relative error of the solution at sampled locations is 
roughly $7\times 10^{-10}$ for both tests.  
Recall, we assume the HBS representation of $\mtx{A}_{oo}$ and its inverse are available. Thus 
the time needed to construct these is not included in our results.  
The results in the first row of Table \ref{tab:precondition_on_refinedfish} are for when the fast matrix vector multiplication for $\tilde{\mtx{A}}$ uses the 
precomputed HBS representation of $\mtx{A}_{oo}$. 
The time for constructing the efficient forward apply of the ELS 
$\left( \tilde{\mtx{A}} +\mtx{LR}\right)$ is 0.53 second, which
includes the construction of $\mtx{A}_{pp}$ and the low-rank factorization $\mtx{Q}\approx\mtx{LR}$.
   As expected the number of iterations is the same as in Table \ref{tab:conditiontest}.  The 
second row \agcmt{in Table \ref{tab:precondition_on_refinedfish}} presents the results when the preconditioner is used.  
The extra time required to construct the preconditioner $T_{\rm pre}$, i.e., for constructing $\mtx{A}^{inv}_{\rm ext}$, includes everything else that was not included in \yzcmt{constructing the efficient forward apply of the ELS 
$\left( \tilde{\mtx{A}} +\mtx{LR}\right)$}
such as the construction and inversion of the Woodbury operator. 
Again the results are comparable to the 
the results in the previous section.  
The preconditioner reduces the number of iterations from 55 to 2, resulting in a 82.7\% reduction in solve time.
And the extra cost for building the preconditioner is justified for problems \agcmt{involving} more than one right-hand-side.

If the problem is not well-conditioned, then the preconditioner may be necessary to obtain an accurate solution with a limited amount of computational resources.
 
\begin{table}[h]
\begin{center}
\begin{tabular}{|l|c|c|c|c|}
\hline
Method & $n_{\rm iter}$ &  $T_{\rm pre}$ & $T_{\rm sol}$ \\
\hline
GMRES with fast mat-vec& 55 & NA & 4.8e-1\\
GMRES with preconditioner & 2& 7.2e-1& 8.3e-2\\
\hline
\end{tabular}
\caption{ Number of iterations $n_{\rm iter}$,
 time in seconds for computing the preconditoner $T_{\rm pre}$ and time in 
seconds for the iterative solver to converge $T_{\rm sol}$ when applying the ELS preconditioner 
to the boundary value problem on the refined fish geometry in Figure \ref{fig:conditiontest}.  The 
red portion of the boundary is refined.  Originally there were $N_c = 128$ points on the red portion.  In the 
new problem there are $N_p = 1024$ points on the red portion of the boundary.  The number of points unchanged
is $N_k=3072$. We assume an HBS representation and the inverse for the original problem are available.
}
\label{tab:precondition_on_refinedfish}
\end{center}
\end{table} 

%

 \input{numerics_revised}

\section{Conclusions}
\label{sec:conclusion}
This manuscript presented a fast direct solver for Stokes BIEs on locally refined discretizations.  \agcmt{The technique makes use of an extended linear system that allows for precomputed fast direct solvers on the unrefined geometry to be utilized.
The numerical results illustrate the new solver's} performance on particulate flow simulations.

\agcmt{For general Stokes problems, two solution approaches are explored.  Which solution 
technique should be used depends on the conditioning of the problem and how many digits
are desired.  For well-conditioned problems, the proposed fast direct solver works extremely well. When the problem has poor conditioning, the fast direct solver
may lose a couple of digits (relative to the compression accuracy).  These digits
can be recovered by using the second solution technique presented here, which is to 
utilize an iterative solver where  
the fast direct solver for the linear system serves as a preconditioner and the compressed representation of the ELS provides the fast matrix
vector multiply.}
Both \agcmt{solution techniques} scale linearly with the size of \agcmt{the unrefined} discretization.
Linear scaling with respect to the number of unknowns added in the local refinement can also be achieved but is not necessary for the considered applications since \agcmt{a relatively low number of points are added}.
%
Numerical examples demonstrated significant speedups; in one test case, the proposed direct solver is roughly 55 times faster than the standard approach. 
For problems with large condition number, more accurate solution may be obtained by using the  proposed preconditoner as compared to the direct solver.  In another test example, the preconditioned GMRES solve for the ELS reduced the number of iterations by a factor of 19 (and total solve time by 3.6X). 
%
%
Our immediate future directions include incorporating close evaluation schemes and extension to three-dimensional problems. 
\section{Acknowledgments}
The authors thank Hai Zhu for providing the Fallopian geometry and Manas Rachh for providing the implementation of the smoothing technique used in the numerical experiments. This work was partially supported by the NSF under grant DMS-2012424.

\bibliographystyle{ieeetr}
\bibliography{main_ref}

\appendix

\section{Proof of Lemma \ref{thm:bound} (Lemma 1 in \cite{YIP1986})}
\begin{proof}

The matrix $ \tilde{\mtx{A}}$ is invertible since it is block diagonal with each block invertible by our assumption.
Let the pseudo-inverses of $\mtx{L}$ and $\mtx{R}$ be defined as above.\\

 It is easy to verify that $\mtx{L}^\dag \mtx{L} = \mtx{R} \mtx{R}^\dag =\mtx{I}$ with dimension $k\times k$.
 By right-multiplying both sides of $\hat{\mtx{A}}_{ext}= \tilde{\mtx{A}}+\mtx{LR}$ by 
 $\tilde{\mtx{A}}^{-1}\mtx{L}$, we get
 $$
 \hat{\mtx{A}}_{ext}\tilde{\mtx{A}}^{-1}\mtx{L}
 = \left(\tilde{\mtx{A}}+\mtx{LR}\right)\tilde{\mtx{A}}^{-1}\mtx{L}
 =\mtx{L}+\mtx{LR}\tilde{\mtx{A}}^{-1}\mtx{L}
 =\mtx{L} \left(\mtx{I} +\mtx{R}\tilde{\mtx{A}}^{-1}\mtx{L}  \right).
 $$
 Now the left-multiplication on both sides of the previous equality by 
 $ \mtx{L}^\dag$ results in the following:
 $$
  \mtx{L}^\dag \hat{\mtx{A}}_{ext}\tilde{\mtx{A}}^{-1}\mtx{L}
  = \mtx{L}^\dag\mtx{L} \left(\mtx{I} +\mtx{R}\tilde{\mtx{A}}^{-1}\mtx{L}  \right)
  =\mtx{I} +\mtx{R}\tilde{\mtx{A}}^{-1}\mtx{L}.
 $$
 Simplifying utilizing the basic properties of pseudo-inverse gives the following
 expression:
$$
 \mtx{L}^\dag \hat{\mtx{A}}_{ext}^{-1}\tilde{\mtx{A}}\mtx{L}
  =\left( \mtx{I} +\mtx{R}\tilde{\mtx{A}}^{-1}\mtx{L} \right)^{-1}.
$$
 Therefore, the condition number for the Woodbury operator is bounded above by
 \begin{equation*}
 \begin{split}
 \kappa\left( \mtx{I} +\mtx{R}\tilde{\mtx{A}}^{-1}\mtx{L} \right) &= \left\| \mtx{I} +\mtx{R}\tilde{\mtx{A}}^{-1}\mtx{L} \right\|
\left\|\left( \mtx{I} +\mtx{R}\tilde{\mtx{A}}^{-1}\mtx{L} \right)^{-1}\right\|\\
& = \left \|   \mtx{L}^\dag \hat{\mtx{A}}_{ext}\tilde{\mtx{A}}^{-1}\mtx{L} \right\|
\left \|  \mtx{L}^\dag \hat{\mtx{A}}_{ext}^{-1}\tilde{\mtx{A}}\mtx{L} \right\|\\
&\leq \left \|   \mtx{L}^\dag\right\|^2  \left \|   \mtx{L}\right\|^2  \left \| \hat{\mtx{A}}_{ext}^{-1}\right\|
\left \|\hat{\mtx{A}}_{ext}\right\| \left\|\tilde{\mtx{A}}^{-1} \right\|\left\|\tilde{\mtx{A}} \right\|\\
&= \left \|   \mtx{L}^\dag\right\|^2  \left \|   \mtx{L}\right\|^2   \kappa\left( \hat{\mtx{A}}_{ext}\right)  \kappa\left(\tilde{\mtx{A}} \right).
\end{split}
 \end{equation*}
 A similar argument using $\mtx{R}^\dag$ gives
  \begin{equation*}
 \begin{split}
 \kappa\left( \mtx{I} +\mtx{R}\tilde{\mtx{A}}^{-1}\mtx{L} \right)
&\leq \left \|   \mtx{R}^\dag\right\|^2  \left \|   \mtx{R}\right\|^2   \kappa\left( \hat{\mtx{A}}_{ext}\right)  \kappa\left(\tilde{\mtx{A}} \right).
\end{split}
 \end{equation*}
 Combining the two bounds above gives equation (\ref{equ:woodburyConditionThm}).
 
\end{proof}

\section{Extended system for the channel with added holes problem}
\label{appendix:variant_ext_sys}
Let $\Gamma_k$ be the original channel bounday given in Figure \ref{fig:ryan_pipe}(a) and $\Gamma_p$ be the union of the holes added in \ref{fig:ryan_pipe}(b). 
Following this subscript notation, the discretized BIE on the ``channel-with-holes" geometry can be reordered into the same format as in (\ref{equ:perturbedlLinearSystem}).
Since no points are deleted, (\ref{equ:perturbedlLinearSystem}) itself serves as an ELS for this problem, and it can be written as
\begin{equation}\label{eq:simple_ext}
\mtx{A}_{nn}=\mtx{A}_{\rm ext} =
\begin{bmatrix}
\mtx{A}_{kk} & \mtx{A}_{kp}\\
\mtx{A}_{pk} &  \mtx{A}_{pp}\\
\end{bmatrix}
=\begin{bmatrix}
\mtx{A}_{kk} & \mtx{0}\\
\mtx{0} &  \mtx{A}_{pp}\\
\end{bmatrix}+\mtx{Q}
\end{equation}
where the update matrix $\mtx{Q}$ can be approximated by a low-rank factorization
$$\mtx{Q}
=\begin{bmatrix}
 \mtx{0}& \mtx{A}_{kp}\\
\mtx{A}_{pk} &   \mtx{0}\\
\end{bmatrix}
\approx \mtx{LR} .$$
If a fast direct solver is already constructed for the original channel geometry, i.e., an approximation to $\mtx{A}_{kk}^{-1}$ is available, then
the solution to (\ref{eq:simple_ext}) can be quickly obtained by a Woodbury formula as described in section \ref{sec:ext}.
The construction of the low-rank approximation for the update matrix is also simpler for this particular problem, since only two subblocks need to be handled.

\end{document}

%% file: numerics_revised.tex
\section{Numerical experiments}
\label{sec:numerics}

This section illustrates the performance of \agcmt{the proposed solution techniques 
for Stokes problems involving locally refined discretizations.}  
The \agcmt{fast direct solver} scales linearly with respect to the number of points in the original discretization and is cheaper than building a fast direct solver from scratch for the new discretization. 
Section \ref{sec:scaling} illustrates the \agcmt{performance of the }fast direct solver when applied to a locally refined channel.  This 
example is from  \cite{RYAN2020ARXIV}. Section \ref{sec:refined_bie} reports on the performance of the fast 
direct solver as a preconditioner when the geometry is complex.  Finally Section \ref{sec:star_lattice} illustrates
the performance of the fast direct solver as a preconditioner when there are a sequence of local refinements 
for the same original geometry.  Such an example arises in many applications including simulations of 
microfluidic devices.

For all test problems, 
the right-hand-side of the BVPs is generated from a known flow and the solution error is 
the average of relative error at chosen target locations in the domain.
All boundaries are discretized via the Nystr\"om method with 16-point composite Gaussian quadrature,
and generalized Gaussian quadrature corrections \cite{quad_review_paper} are used to handle the weakly singular kernels.
The solver also works with other quadrature corrections, such as \cite{Alpert99, helsing_corner,kapur_rokhlin}.

All experiments were run on a dual 2.3 GHz Intel Xeon Processor E5-2695 v3 desktop workstation with 256 GB of RAM. The code is implemented in MATLAB, apart from the interpolatory decomposition routine, which is in FORTRAN.

\yzcmt{
To illustrate the performance of the solver, 
we introduce the following notations for reporting times and errors.
For notation consistency, we use regular capital letters such as $T$ and $E$ for problems defined
on the original discretization (or geometry) and letters with tilde, such as $\tilde{T}$
and $\tilde{E}$ for problems on the locally refined discretization (or perturbed geometry).
For the problem on the original discretization (or geometry), we define} 
\yzcmt{
\begin{itemize}
    \item $T_{\rm HBS,\,comp}$ and $T_{\rm HBS,\,inv}$:  the time in seconds for building the HBS compression of the discretized boundary integral operator and that for inverting the compression, i.e., building the HBS inverse, respectively.
    \item $T_{\rm HBS,\,Dsol}$: the time in seconds for applying the HBS inverse to a given right-hand-side vector. ``Dsol" stands for ``one direct solve".
    \item $T_{\rm HBS,\,Gsol}$: the time in seconds for solving for one right-hand-side vector using GMRES with HBS compression accelerated matrix-vector multiplication. ``Gsol" stands for ``one GMRES solve".
    \item $T_{\rm HBS,\,PGsol}$: the time in seconds for solving for one right-hand-side vector using a preconditioned GMRES with HBS compression accelerated  matrix-vector multiplication, where the HBS inverse is used as the preconditioner.   ``PGsol" stands for ``one preconditioned GMRES solve".
    \item $E_{\rm HBS,\, Dsol}$, $E_{\rm HBS,\, Gsol}$ and $E_{\rm HBS,\, PGsol}$: the average relative error at sample domain locations for the three different solve options respectively.
\end{itemize}}
\yzcmt{
For the problem on the locally refined discretization (or perturbed geometry), we define
\begin{itemize}
    \item $\tilde{T}_{\rm HBS,\,comp}$, $\tilde{T}_{\rm HBS,\,inv}$, $\tilde{T}_{\rm HBS,\,Dsol}$, $\tilde{T}_{\rm HBS,\,Gsol}$, and $\tilde{T}_{\rm HBS,\,PGsol}$: time in seconds  similar to those categories for the original discretization (or geometry).
    \item $\tilde{E}_{\rm HBS,\, Dsol}$, $\tilde{E}_{\rm HBS,\, Gsol}$ and $\tilde{E}_{\rm HBS,\, PGsol}$: error similar to those categories for the original discretization (or geometry).
    \item $\tilde{T}_{\rm ELS,\,comp}$: the time in seconds for building $\mtx{A}_{pp}$ and $\mtx{LR}\approx \mtx{Q}$ in formulating the fast ELS approximation. Note we assume a HBS compression for $\mtx{A}_{oo}$ is available. 
    \item $\tilde{T}_{\rm ELS,\,inv}$: the time in seconds for building the operators needed in the Woodbury formula for applying the inverse approximation of the ELS: $\mtx{A}_{pp}^{-1}$, $\tilde{\mtx{A}}^{-1}\mtx{L}$, $\mtx{R}\tilde{\mtx{A}}^{-1}\mtx{L}$ and $\left(\mtx{I}+\mtx{R}\tilde{\mtx{A}}^{-1}\mtx{L}\right)^{-1}$. Note we assume a HBS inverse approximation for $\mtx{A}_{oo}$ is available. 
    \item $\tilde{T}_{\rm ELS,\,Dsol}$: the time in seconds for applying the approximate ELS inverse $(\tilde{\mtx{A}}+\mtx{LR})^{-1}$ via the Woodbury formula to a given right-hand-side vector.
    \item $\tilde{T}_{\rm ELS,\,Gsol}$: the time in seconds for solving the approximate ELS $(\tilde{\mtx{A}}+\mtx{LR})\bm{\tau}_{\rm ext}=\vct{g}_{\rm ext}$ for one right-hand-side vector $\vct{g}_{\rm ext}$ using GMRES.
    \item $\tilde{T}_{\rm ELS,\,PGsol}$: the time in seconds for solving the approximate ELS $(\tilde{\mtx{A}}+\mtx{LR})\bm{\tau}_{\rm ext}=\vct{g}_{\rm ext}$ for one right-hand-side vector $\vct{g}_{\rm ext}$ using GMRES, where the approximate ELS inverse is used as the preconditioner. 
    \item $\tilde{E}_{\rm ELS,\, Dsol}$, $\tilde{E}_{\rm ELS,\, Gsol}$ and $\tilde{E}_{\rm ELS,\, PGsol}$: the average relative error at sample domain locations for the three different ELS solve options respectively.
\end{itemize}}
\yzcmt{
The accuracy for HBS compression and low-rank approximation is set to $10^{-10}$ unless specified otherwise.
}


\subsection{Asymptotic scaling experiments}
\label{sec:scaling}
\begin{figure}[ht]
\begin{center}
\begin{tabular}{c c c}
\includegraphics[height=45mm,trim={3cm 0 3cm 0}, clip]{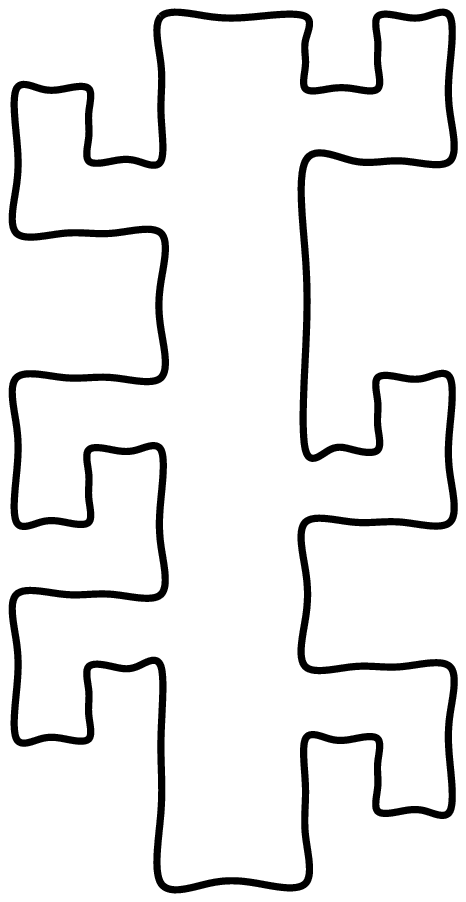}
 &
\includegraphics[height=45mm,trim={3cm 0 3cm 0}, clip]{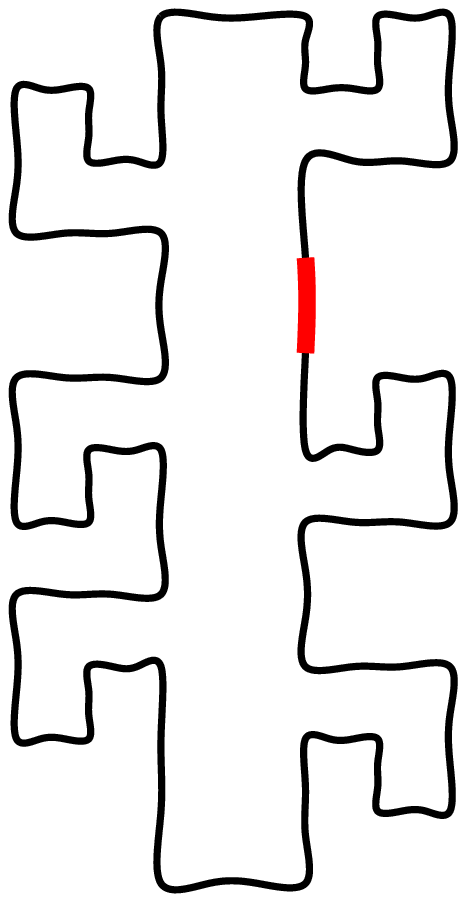}
 &
\includegraphics[height=45mm,trim={3cm 0 3cm 0}, clip]{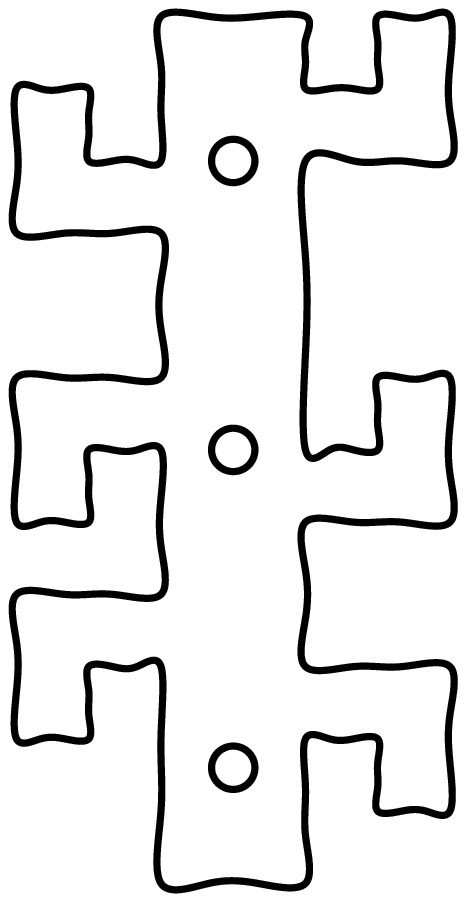}
 \\
(a)&(b)&(c)\\
\end{tabular}
\end{center}
\caption{ (a) The original channel geometry. (b) The channel geometry with a locally refined segment highlighted in red.
(c) The channel geometry with three interior holes added.  }
\label{fig:ryan_pipe}
\end{figure}

\agcmt{This section illustrates the performance of the fast direct solver for the ELS 
(presented in Section \ref{sec:fdsolver}) when applied to a} Stokes problem with a confined geometry with two types of
modifications: locally refining a part of the boundary and adding holes.  
Figure \ref{fig:ryan_pipe}(a) 
illustrates the channel. Figure \ref{fig:ryan_pipe}(b) and (c) 
illustrate the modification of local refinement (in red) and adding holes, respectively.  
The geometry is generated by applying cubic splines with periodic conditions to 121 spline knot locations 
 (with the first and last knots give the same physical point on the geometry) and was first seen in \cite{RYAN2020ARXIV}.
 The channel is discretized by \agcmt{using the same number} of Gaussian panels per subinterval in \agcmt{in the} cubic spline \agcmt{geometry generation.}
For example, the total number of discretization points on the channel $N_{\rm channel}=1920$ corresponds to $120$ Gaussian panels in total and 1 panel per subinterval.  \agcmt{If there are two panels per subinterval, the number of discretization points doubles.}
The circular holes are each discretized with 10 panels \agcmt{which means there are} $160$ quadrature points \agcmt{per circle}. 

\begin{remark}
 The addition of holes is similar to the original examples used in \cite{RYAN2020ARXIV} and fits in the definition of \agcmt{a} locally perturbed geometry 
 as defined in \cite{ZHANG2018, ZHANG2021}. However, the extended system is slightly different from the one given in section \ref{sec:ext}
 as we are only adding points \agcmt{for the new 
 boundary and there is no deletion or cutting of points on the original geometry}. The corresponding ESL formulation is given in Appendix \ref{appendix:variant_ext_sys}.
\end{remark}

 The Dirichlet boundary data for the interior channel BVP (Figure \ref{fig:ryan_pipe}(a) and (b)) 
is generated by 5 exterior Stokeslets outside of the channel geometry. For these two problem, the 
solution is represented with the double layer kernel (as discussed in Section \ref{sec:int}).  
The Dirichlet boundary data for the BVP with holes (Figure \ref{fig:ryan_pipe}(c))
is generated by the same 5 exterior Stokeslets outside the channel geometry and \yzcmt{five
additional Stokeslets placed inside the added holes (two stokelets per hole for the bottom two holes and one stokelet in the top hole).} 
The solution interior to the channel and exterior to the holes is represented as a double layer potential on the channel plus a combined field potential on the holes.

The observed condition number \agcmt{of} the discretized integral operator for all the problems in this section is on the order of $10^5$.
The condition number of the Woodbury operator is on the order of $10^3$ for the problem with the added holes and $10^5$ 
for the problem with the local refinement.
\yzcmt{The observed rank numbers for the low-rank approximation of the update matrix $\mtx{Q}$,
which is also the size of the Woodbury system, is roughly  60 for the problem with the local refinement and 340 for the problem with the added holes.}

\yzcmt{
Let \agcmt{$N_{\rm channel}$ denote }the number of discretization points on the original channel.
Table \ref{tab:ryan_pipe_original} \agcmt{reports on} the performance of the HBS solver applied to the original geometry (and discretization).
For the locally refined discretization, let \agcmt{$N_c$ and $N_p$ denote the number of points removed and added, respectively}.
For the channel with holes geometry, let $N_{\rm holes}$ be the total number of discretization points placed on the three holes.}
\yzcmt{
 The results from Table \ref{tab:ryan_pipe_refine} and \ref{tab:ryan_pipe_holes} \agcmt{report on} the performance of the 
proposed ELS formulation based fast direct solver applied to the geometries in Figures \ref{fig:ryan_pipe}(b) and (c).
The size of each test cases is given by the total degree of freedom, which is double the number of 
discretization points.
To show the scaling of the ELS fast direct solver,
the values for $N_{\rm channel}$, $N_c$ and $N_p$ are all doubled as the test size increases.}
\yzcmt{
Both the HBS solver and the proposed fast direct solver scale linearly with respect to the number of points on the channel 
geometry.  The cost of using Algorithm 1 is significantly less than building the original HBS solver.  
This \agcmt{means} that Algorithm 1 is more computationally efficient than building a fast direct solver \agcmt{from scratch} for the new discretization.  \agcmt{It is worth noting that}  the time required for building the
ELS compression for the addition of holes example is much higher than that for refining the channel boundary given the same $N_{\rm channel}$. 
For example, when $2N_{\rm channel}=122880$, ELS compression for adding holes is about 8 times of that for refining a segment, although the points added for the holes $N_{\rm holes}$ is only 1/3 of the points added $N_p$ due to the refinement.}
\yzcmt{
 This is due to the fact that change to the system for adding the three holes is ``less local'' than that for refining a segment of the channel, resulting in much higher rank numbers and more expensive compression of the update matrix $\mtx{Q}$.  
 \agcmt{For} the same reason, the time required for applying the inverse of the ELS 
 \agcmt{when} adding holes is also more than that for refining a segment of the channel.}

\yzcmt{
The solution error for all test cases \agcmt{is} maintained at $10^{-10}$ \agcmt{since} the geometry is fully resolved and the tolerance for HBS compression and low-rank approximations is set to be $10^{-10}$.
}


\begin{table}[h]
\begin{center}
\begin{tabular}{|c||c|c|c||c|}
\hline
&&&&
\\[-1em]
2$N_{\rm channel}$ & $T_{\rm HBS,\, comp}$ & $T_{\rm HBS,\, inv}$ &$T_{\rm HBS,\, sol}$ & $E_{\rm HBS,\, Dsol}$\\
\hline
30720  &65.7 & 7.0& 0.070&  1.1e-10\\
61440  &90.7 &9.9 &0.140 & 1.4e-10 \\
122880  &132.8 &15.8 & 0.264& 3.22e-10   \\
\hline
\end{tabular}
\end{center}
\caption{ The time in seconds and error for using HBS compression and inversion to solve the BIE on the channel geometry  with the original discretization (illustrated 
in Figure \ref{fig:ryan_pipe}(a)). 
 The number of discretization points on the channel is $N_{\rm channel}$. 
The size of the linear system is $2N_{\rm channel}\times 2N_{\rm channel}$.
}
\label{tab:ryan_pipe_original}
\end{table}

\begin{table}[h]
\begin{center}
\begin{tabular}{|c||c|c|c||c|}
\hline
&&&&
\\[-1em]
2$N_{\rm channel}$,  2$N_c$,  2$N_p$ &  $\tilde{T}_{\rm ELS,\, comp}$& $\tilde{T}_{\rm ELS,\, inv}$ & $\tilde{T}_{\rm ELS,\, Dsol}$  & $\tilde{E}_{\rm ELS,\, Dsol}$\\
\hline
30720, 192, 768  &1.4 &0.8 & 0.088& 2.5e-10 \\
61440, 384, 1536 & 3.0&1.2 &0.113 &  5.8e-10 \\
122880, 768, 3072   &7.1 & 2.3  &0.184 &  4.8e-10 \\
\hline
\end{tabular}
\end{center}
\caption{ The time in seconds and error for using the Woodbury formula to solve the ELS for the boundary value problem on
the channel with local refinement illustrated in Figure \ref{fig:ryan_pipe}(b).
  The
number of discretization points cut and added on the red portion of the boundary are $N_c$ and $N_p$ respectively.
}
\label{tab:ryan_pipe_refine}
\end{table}

\begin{table}[h]
\begin{center}
\begin{tabular}{|c||c|c|c||c|}
\hline
&&&&
\\[-1em]
2$N_{\rm channel}$, 2$N_{\rm holes}$ &  $\tilde{T}_{\rm ELS,\, comp}$& $\tilde{T}_{\rm ELS,\, inv}$ & $\tilde{T}_{\rm ELS,\, Dsol}$  & $\tilde{E}_{\rm ELS,\, Dsol}$\\
\hline
30720, 960 &14.9 &2.1 & 0.059&9.8e-11\\
61440, 960  & 28.8&4.3 &0.103 &2.3e-10\\
122880, 960  &56.4 & 7.9&0.296 &2.1e-10\\
\hline
\end{tabular}
\end{center}
\caption{The time in seconds and error for using the Woodbury formula to solve the ELS for the boundary value problem on
the channel with three interior holes illustrated in Figure \ref{fig:ryan_pipe}(c). 
The circular holes in Figure \ref{fig:ryan_pipe}(b) are each discretized with 10 panels and $160$ quadrature points, resulting a total of $N_{\rm holes}=480$ points.  }
\label{tab:ryan_pipe_holes}
\end{table}


\subsection{Complex geometry with local refinement}
\label{sec:refined_bie}
\begin{figure}[ht]
\begin{center}
\includegraphics[height=80mm,trim={1cm 1cm 1cm 1cm}, clip]{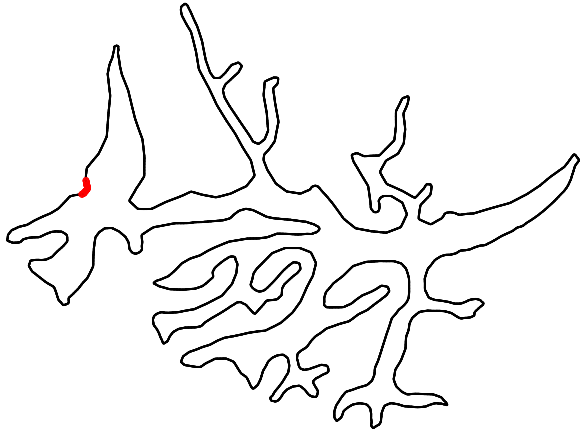}
\end{center}
\caption{ Partial Fallopian tube based on data extracted from the experiments in \cite{Guo_phyrev2020}.
A small segment highlighted in red is chosen to be locally refined. The geometry is generated by \cite{Beylkin_curve_fitting}.}
\label{fig:Fallopian}
\end{figure}

This section considers an interior problem on the complex Fallopian tube geometry illustrated in Figure \ref{fig:Fallopian}. 
The geometry is created by extracting data points from \yzcmt{Figure 1 of} 
\cite{Guo_phyrev2020} and connecting them smoothly via the technique in \cite{Beylkin_curve_fitting}.  
The solution to the problem is generated by \agcmt{placing} Stokeslets on the exterior of the geometry. \agcmt{The boundary data is generated via
this known solution.}
\yzcmt{Discretizing \agcmt{the} complex geometry \agcmt{in Figure  \ref{fig:Fallopian}} 
results in an integral equation with \agcmt{a} high condition number.
\agcmt{An iterative solver requires a large number of iterations in order to converge. } 
}
\yzcmt{\agcmt{The experiments in this section discretize the original Fallopian tube boundary (pre-refinement) } with 1600 Gaussian panels (25600 points and 51200 degrees of freedom),
which results in relative error \agcmt{of} approximately $4\times 10^{-5}$.
To understand the conditioning of the linear system, we consider the smallest matrices 
that are inverted in the hierarchical tree using the HBS solver.  These matrices (corresponding
to the first three levels in the tree) have condition numbers on the order of 
$10^{8}$ to $ 10^{11}$. }

\yzcmt{
For the refined discretization problem, the red portion of the boundary highlighted 
in Figure \ref{fig:Fallopian} goes from having $6$ panels to $24$.  Since a 16 point 
Gaussian quadrature is used, the number of points kept, cut and added are $N_k=25,504$, $N_c=96$, and $N_p=384$, respectively.  
The iterative solver stops when the 
relative residual is on the order of $10^{-6}$.  
For the boundary integral equation on the original discretization, we either build only HBS representation of the discretized boundary integral equation and couple it with GMRES
or also build the HBS inverse and apply it directly to the given right-hand-side.
For the refined problem, we consider the discretized BIE and the equivalent ELS and a fast direct solver and an iterative solver for each.
Additionally, we also use the direct solver for the ELS, built as Algorithm 1, to precondition the GMRES solve.
}

\yzcmt{
Table \ref{tab:Fallopian_Dsol} reports the time required to solve the BIE on the original discretization,
the BIE on the refined discretization, and the approximate ELS on the refined discretization using a fast direct solver. 
The total time for precomputation includes two parts: the forward compression indicated by subscript notation ``comp" and the inversion indicated by the subscript notation ``inv". 
\agcmt{Tables \ref{tab:Fallopian_Dsol}(b) and  \ref{tab:Fallopian_Dsol}(c) demonstrate that for
this geometry the proposed direct solver for the ELS is more efficient than building a HBS solver from scratch for the refined problem. }}
\yzcmt{
\agcmt{In fact}, the cost for constructing a forward compression for the ELS for the refined problem is only 1.3\% of the cost of constructing a HBS from scratch.
The cost of constructing the inverse operator is only 7.3\% of that of HBS inverse. 
}

\yzcmt{
Table \ref{tab:Fallopian_Gsol} reports the time in seconds for 
 the unpreconditioned GMRES approach for the original and refined problems.   The precomputation for this approach only involves \agcmt{the compression of the forward operator} and is lower than that for the direct solution approach \agcmt{since an approximate inverse is not constructed.}}
\yzcmt{
 However, due to the \agcmt{poor} conditioning of the problem, more than 500 GMRES iterations are required to reach the desired tolerance of $10^{-6}$.  \agcmt{The time required to solve the integral equation via unprecondtioned GMRES is much higher
 for each new right-hand-side than the direct solver. 
 Tables \ref{tab:Fallopian_Gsol}(b) and (c) show that 
  solving the  
  approximate} ELS (\ref{eq:ELS_approximate}) is two orders of magnitude cheaper than building the HBS compression of the BIE for the refined problem. 
  Applying the forward operator for the ELS is slightly more expensive than applying the HBS forward compression. 
 }
 
 \yzcmt{
 Table \ref{tab:Fallopian_PGsol} reports the time in seconds for the preconditioned GMRES approach applied to the approximate ELS. Here Algorithm 1, i.e., the inverse of the ELS obtained by the Woodbury formula, is used to precondition the fast representation of the ELS. 
 The precomputation time of this approach is equal to that of the direct \agcmt{solver} approach for the ELS. 
 The number of GMRES iterations required \agcmt{for the convergence criterion to be met is} reduced from 520 to 6, leading to a significant reduction in total cost
 even for only one right-hand-solve when compared to the results in Table \ref{tab:Fallopian_Gsol}(c).
The cost of solving one additional right-hand-side vector via the preconditioned GMRES approach for the ELS is about 5.4\% of that via unpreconditioned GMRES approach. 
}

\begin{table}[h]
  \begin{center}
 \begin{tabular}{|c|c|c|}
  \hline
  $T_{\rm HBS, \, comp}$ & $T_{\rm HBS, \, inv}$ &$T_{\rm HBS, \, Dsol}$\\
  \hline
  4.09e+2 & 2.76e+1 & 7.87e-2\\
    \hline
  \end{tabular}
   
  (a) 
   
  \vspace{1em}
   
  \begin{tabular}{|c|c|c|}
  \hline
  &&
  \\[-1em]
  $\tilde{T}_{\rm HBS, \, comp}$ & $\tilde{T}_{\rm HBS, \, inv}$ &$\tilde{T}_{\rm HBS, \, Dsol}$\\
  \hline
  4.09e+2 & 2.76e+1 & 7.87e-2\\
    \hline
  \end{tabular}
   
  (b)
   
  \vspace{1em}
   
  \begin{tabular}{|c|c|c|}
  \hline
  &&
  \\[-1em]
  $\tilde{T}_{\rm ELS, \, comp}$ & $\tilde{T}_{\rm ELS, \, inv}$ &$\tilde{T}_{\rm ELS, \, Dsol}$\\
  \hline
  5.33e+0 & 2.01e+0 & 9.62e-2\\
    \hline
  \end{tabular}
   
  (c)
   
  \caption{Time in seconds for solving (a) the original and (b) the refined problem defined on the  Fallopian tube geometry (Figure \ref{fig:Fallopian}) via HBS inversion. (c) corresponds to the proposed fast direct solver for the ELS of the refined problem. }
    \label{tab:Fallopian_Dsol}
    \end{center}
  \end{table}
  \begin{table}[h]
  \begin{center}
 \begin{tabular}{|c|c|}
  \hline
  $T_{\rm HBS, \, comp}$ &$T_{\rm HBS, \, Gsol}$\\
  \hline
  4.09e+2 & 5.95e+1 ($n_{\rm iter}=519$)\\
    \hline
  \end{tabular}
   
  (a) 
   
  \vspace{1em}
   
  \begin{tabular}{|c|c|}
  \hline
  &
  \\[-1em]
  $\tilde{T}_{\rm HBS, \, comp}$ &$\tilde{T}_{\rm HBS, \, Gsol}$\\
  \hline
  4.09e+2 &  6.59e+1 ($n_{\rm iter}=519$)\\
    \hline
  \end{tabular}
   
  (b)
   
  \vspace{1em}
   
  \begin{tabular}{|c|c|}
  \hline
  &
  \\[-1em]
  $\tilde{T}_{\rm ELS, \, comp}$ &$\tilde{T}_{\rm ELS, \, Gsol}$\\
  \hline
  5.33e+0 &  7.07e+1 ($n_{\rm iter}=520$)\\
    \hline
  \end{tabular}
   
  (c)
   
  \caption{Time in seconds for solving (a) the original and (b) the refined problem defined on the  Fallopian tube geometry (Figure \ref{fig:Fallopian}) via GMRES with HBS compression accelerated matrix-vector multiplication.
  (c) corresponds to solving the ELS of the refined problem via GMRES. The number of GMRES iterations $n_{\rm iter}$ required to converge to tolerance $10^{-6}$ is also reported.}
    \label{tab:Fallopian_Gsol}
    \end{center}
  \end{table}

  \begin{table}[h]
  \begin{center}
 
  \begin{tabular}{|c|c|c|}
  \hline
  &&
  \\[-1em]
  $\tilde{T}_{\rm ELS, \, comp}$ &$\tilde{T}_{\rm ELS, \, inv}$ &$\tilde{T}_{\rm ELS, \, PGsol}$\\
  \hline
    5.33e+0 &  2.01e+0 &3.80e+0 ($n_{\rm iter}=6$)\\
    \hline
  \end{tabular}

  \caption{Time in seconds for solving the ELS of the refined problem defined on the  Fallopian tube geometry (Figure \ref{fig:Fallopian}) via preconditioned GMRES, where Algorithm 1 is used as the preconditioner. The number of GMRES iterations $n_{\rm iter}$ required to converge to tolerance $10^{-6}$ is also reported.}
    \label{tab:Fallopian_PGsol}
    \end{center}
  \end{table}

\subsection{Relocating region of local refinement}
\label{sec:star_lattice}
\begin{figure}[ht]
\begin{center}
\begin{tabular}{c}
\includegraphics[height=120mm,trim={5cm 2cm 5cm 0}, clip]{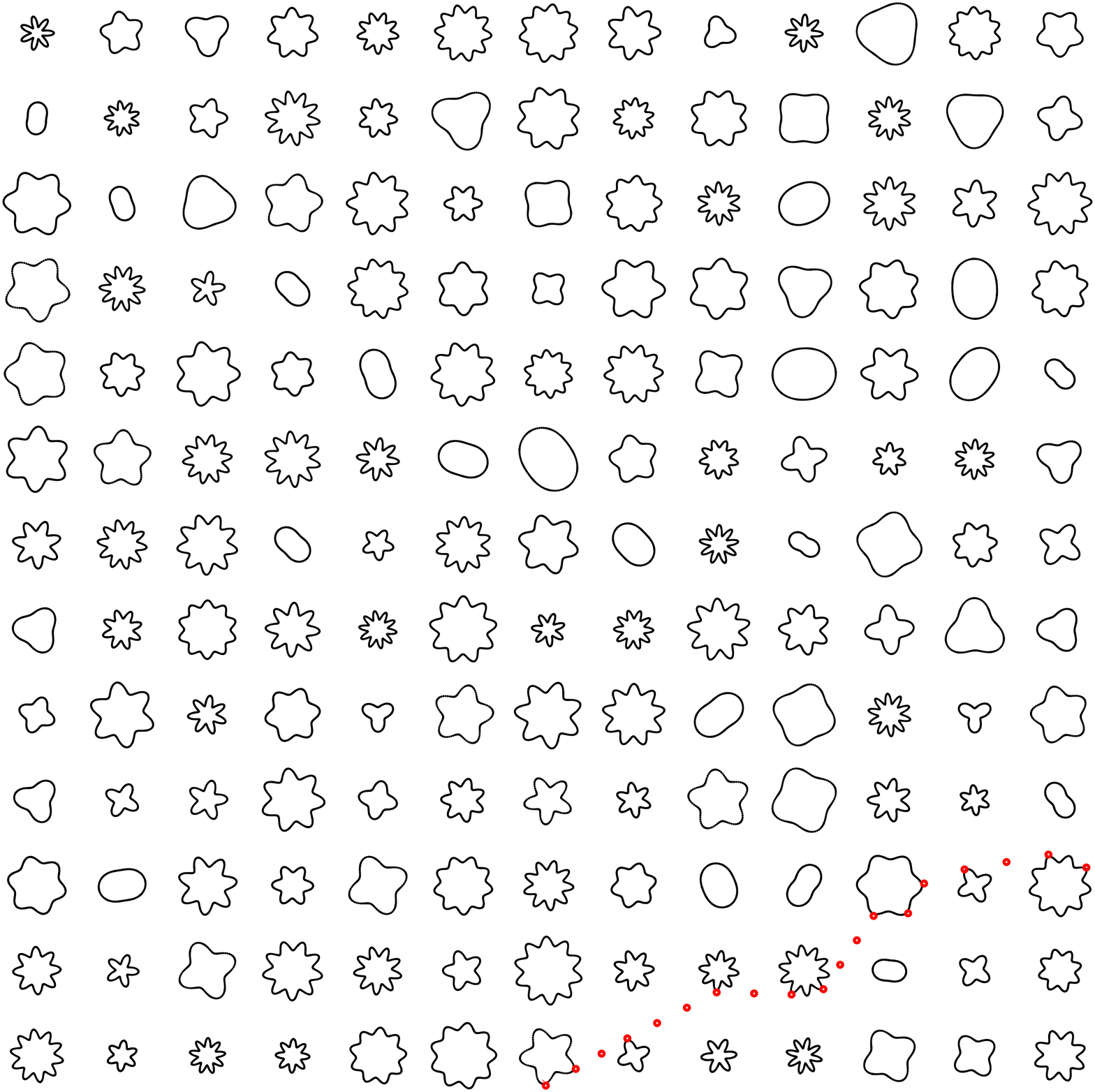}
\\
(a)\\
\\
\includegraphics[height=80mm,trim={0cm 10cm 0cm 6cm}, clip]{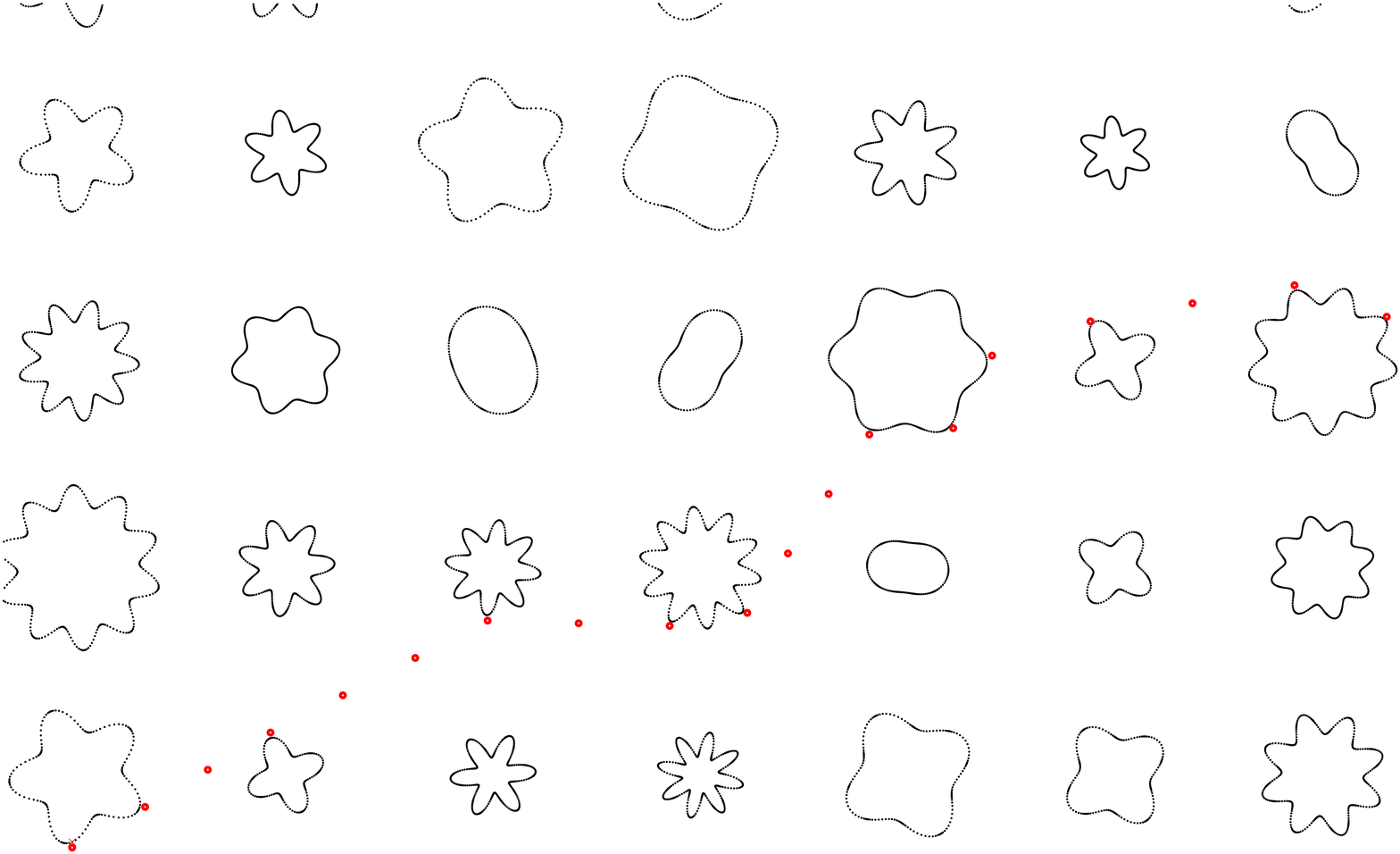}
 \\
(b)\\
\end{tabular}
\caption{ (a) A collection of star-shape obstacles with different snapshots of body locations. (b) Zoomed-in in the region near the snapshots of the body locations. 
The locations are chosen artificially and do not represent any physical movement of body in Stokes flow. 19 locations are chosen, out of which 12 are close to certain part of the obstacle boundary and incurs local refinement of the obstacle discretization.}
\label{fig:star_lattice}
\end{center}
\end{figure}

This section illustrates the potential of using the fast direct solver
presented in Algorithm 1 as preconditioner for many Stokes 
problems involving a body moving through a collection of star-shaped
 obstacles shown in Figure \ref{fig:star_lattice}.  This example is 
 representative of applications \agcmt{such as sorting with} a microfluidic device.
For the original discretization, 10 panels are placed on each star with less or equal to 5 prongs \agcmt{and 20 panels} are placed on stars with more than 5 prongs.  \agcmt{With the 16-point Gaussian quadrature,} this results in a total of 42400 discretization points and \yzcmt{84800 degrees of freedom}.
For demonstration purposes, we do not simulate the true physics of any body moving in the domain; instead, we assume the body appears at certain locations at some time step, \agcmt{as illustrated in} in Figure \ref{fig:star_lattice}. These can be viewed as snapshots of a body moving through the obstacles.
The body moving through the obstacles is much smaller in scale than any of the stars.  Thus the discretization of one or more obstacles will need to be locally refined as the body approaches 
those obstacles. Since the body is moving, the regions of local refinement are expected
to be different for each snapshot.
\yzcmt{Previously refined regions may be coarsened back into the original discretization as the body moves away.}
In this example, 19 snapshot locations are chosen.  In 12 of these snapshots, the 
body is close to an obstacle and local refinement is needed.  In the other 7 snapshots,
no local refinement is needed.  \agcmt{We} consider $5$ different ways
of solving the linear system for the 19 different boundary value problems.  These
solution techniques are: 
\begin{enumerate}[(1)]
\item \texttt{GMRES-indy}: Treat each of the 13 different discretization as independent 
boundary value problems, building a forward HBS representation for each, and using this to accelerate the GMRES solve for each snapshot;
\item \texttt{Direct-indy}: Treat each of the 13 different discretization as independent boundary 
value problems and build a HBS solver for each one; 
\item \texttt{GMRES-Local}: Build a HBS forward representation for the original discretization and use it to accelerate the GMRES solve for the ELS for each problem requiring local refinement;
\item \texttt{Direct-Local}: Build a HBS solver for the original discretization and use it to build a fast direct solver for the ELS according to Algorithm 1 for each problem needing local refinement;
\item \texttt{PGMRES-Local}: Build a HBS solver for the original discretization and use it to precondition the GMRES solve for the ELS for each problem requiring local refinement.
\end{enumerate}

The tolerance for GMRES is set to $10^{-11}$. 
\yzcmt{For the boundary value problems that do not require local refinement, using the HBS matrix-vector 
acceleration of GMRES 
results in a relative error on the order of $10^{-9}$. Using the 
HBS solver loses two digits; i.e., the relative error that results from this solver
is on the order of $10^{-7}$.  
}
Thus for the two techniques (2) and (4) \agcmt{where} the direct solver is used as an actual solver and not a preconditioner, the accuracy is approximately $10^{-7}$. \agcmt{When} the HBS solver or the ELS \agcmt{fast direct} solver in Algorithm 1 is used
\agcmt{as} the preconditioner, the error is approximately $10^{-9}$. 

\yzcmt{To compare efficiency of the five approaches, 
 we first \agcmt{report} the time in seconds for solving the problem on the original discretization and that on one particular refined discretization, which corresponds to the 
first snapshot with the body located at the very bottom left of Figure \ref{fig:star_lattice} (b). The results are presented in Table \ref{tab:star_lattice_Dsol}, \ref{tab:star_lattice_Gsol} and \ref{tab:star_lattice_PGsol} in the same format as 
\agcmt{the corresponding results} for the Fallopian tube geometry in the previous section.
}
Since the direct solver does not achieve the 
full possible accuracy of the discretization, using it as preconditioner is reasonable
and it greatly decreases the number of iterations needed for an iterative solver to
converge.  

\begin{table}[h]
  \begin{center}
 \begin{tabular}{|c|c|c|}
  \hline
  $T_{\rm HBS, \, comp}$ & $T_{\rm HBS, \, inv}$ &$T_{\rm HBS, \, Dsol}$\\
  \hline
  8.23e+2 & 3.38e+1 & 2.57e-1\\
    \hline
  \end{tabular}
   
  (a) 
   
  \vspace{1em}
   
  \begin{tabular}{|c|c|c|}
  \hline
  &&
  \\[-1em]
  $\tilde{T}_{\rm HBS, \, comp}$ & $\tilde{T}_{\rm HBS, \, inv}$ &$\tilde{T}_{\rm HBS, \, Dsol}$\\
  \hline
  8.07e+2 & 3.38e+1 & 2.98e-1\\
    \hline
  \end{tabular}
   
  (b)
   
  \vspace{1em}
   
  \begin{tabular}{|c|c|c|}
  \hline
  &&
  \\[-1em]
  $\tilde{T}_{\rm ELS, \, comp}$ & $\tilde{T}_{\rm ELS, \, inv}$ &$\tilde{T}_{\rm ELS, \, Dsol}$\\
  \hline
  8.39e+0 & 6.43e+0 & 4.18e-1\\
    \hline
  \end{tabular}
   
  (c)
   
  \caption{Time in seconds for solving (a) the original and (b) the refined problem defined on the star-shape obstacle geometry in Figure \ref{fig:star_lattice}  via HBS inversion. (c) corresponds to the proposed fast direct solver for the ELS of the refined problem. }
    \label{tab:star_lattice_Dsol}
    \end{center}
  \end{table}
  \begin{table}[h]
  \begin{center}
 \begin{tabular}{|c|c|}
  \hline
  $T_{\rm HBS, \, comp}$ &$T_{\rm HBS, \, Gsol}$\\
  \hline
  8.23e+2 & 3.27e+1 ($n_{\rm iter}=113$)\\
    \hline
  \end{tabular}
   
  (a) 
   
  \vspace{1em}
   
  \begin{tabular}{|c|c|}
  \hline
  &
  \\[-1em]
  $\tilde{T}_{\rm HBS, \, comp}$ &$\tilde{T}_{\rm HBS, \, Gsol}$\\
  \hline
  8.07e+2 &  3.30e+1 ($n_{\rm iter}=113$)\\
    \hline
  \end{tabular}
   
  (b)
   
  \vspace{1em}
   
  \begin{tabular}{|c|c|}
  \hline
  &
  \\[-1em]
  $\tilde{T}_{\rm ELS, \, comp}$ &$\tilde{T}_{\rm ELS, \, Gsol}$\\
  \hline
  8.39e+0 &  3.46e+1 ($n_{\rm iter}=113$)\\
    \hline
  \end{tabular}
   
  (c)
   
  \caption{Time in seconds for solving (a) the original and (b) the refined problem defined on the star-shape obstacle geometry in Figure \ref{fig:star_lattice} via GMRES with HBS compression accelerated matrix-vector multiplication.
  (c) corresponds to the solving the ELS of the refined problem via GMRES. The number of GMRES iterations $n_{\rm iter}$ required to converge to tolerance $10^{-11}$ is also reported.}
    \label{tab:star_lattice_Gsol}
    \end{center}
  \end{table}

  \begin{table}[h]
  \begin{center}
 
  \begin{tabular}{|c|c|c|}
  \hline
  &&
  \\[-1em]
  $\tilde{T}_{\rm ELS, \, comp}$ &$\tilde{T}_{\rm ELS, \, inv}$ &$\tilde{T}_{\rm ELS, \, PGsol}$\\
  \hline
    8.39e+0 &  6.43e+0 &9.50e+0 ($n_{\rm iter}=6$)\\
    \hline
  \end{tabular}

  \caption{Time in seconds for solving the ELS of the refined problem defined on the star-shape obstacle geometry in Figure \ref{fig:star_lattice} via preconditioned GMRES, where Algorithm 1 is used as the preconditioner. The number of GMRES iterations $n_{\rm iter}$ required to converge to tolerance $10^{-11}$ is also reported.}
    \label{tab:star_lattice_PGsol}
    \end{center}
  \end{table}

\yzcmt{With the step-by-step cost summarized in Table \ref{tab:star_lattice_Dsol}, \ref{tab:star_lattice_Gsol} and \ref{tab:star_lattice_PGsol}, 
we can approximate the total cost for each of the five approaches handling all 19 snapshots by simple addition and multiplication,
assuming that the cost for solving the ELS for each snapshot that requires a refinement is the same.}
To get an idea of the speed up for solving problems involving the 19 multiple snapshots given in Figure  \ref{fig:star_lattice}(a), Table \ref{tab:star_lattice2} collects the time necessary 
for each part of the $5$ solution techniques.  
The different times reported are:
\begin{itemize}
\item $T_{\rm static}$: The time in seconds for constructing any of the operators 
needed for the solution technique on the original discretization.  For techniques (1) 
and (3), only constructing an approximation of $\mtx{A}_{oo}$ via HBS is needed. For 
the other options, the construction of the approximate inverse of $\mtx{A}_{oo}$ is 
also needed.  This is a ``static'' computation since it is independent of future time steps and potential local refinement.
\item $T_{\rm Osol}$: The time in seconds for solving a problem where local refinement is 
not needed.  ``Osol" stands for ``solve for the original discretization"
\item $T_{\rm Rsol}$: The time in seconds for solving a problem where local refinement is 
needed. ``Rsol" stands for ``solve for one refined discretization". 
\end{itemize}

Approaches (3-5) which utilize the ELS are more efficient than 
building new HBS solver from scratch each time or only when there is local refinement. %
For these experiments Approach (4) is the most efficient but if the fully attainable accuracy
is desired, Approach (5) should be used as it is both efficient and accurate. 
The \agcmt{previous standard solution technique for this type of problem was Approach (1).}  The 
proposed direct solver (4) and the proposed preconditioned solver (5) are 127 and 3.5 times
faster than Approach (1) when local refinement is not needed.  When refinement 
is needed, Approaches (4) and (5) are 55 and 34.6 times faster than Approach (1), respectively.
Since the applications of interest (such as \cite{MARPLE2016_periodicstokes}) involve hundreds to thousands of solves, it is
definitely worth using the ELS.  If the user is okay losing
a couple of digits, the fast direct solver is an ideal choice.  If the digits are needed,
then the preconditioned iterative solver is still going to be significantly faster
than Approach (1).

\begin{remark}
 The dominate cost $T_{\rm Rsol}$ for the ELS solution techniques 
 is the cost of creating the low rank factorization of $\mtx{Q}$.  In most applications,
 several snapshots can use the same refinement and thus the same factorization of 
 $\mtx{Q}$.  The reuse of the factorization will decrease $T_{\rm Rsol}$ significantly.
For example, \agcmt{in the experiments corresponding to} the two body locations on the left bottom of Figure \ref{fig:star_lattice},
two different regions of the same five-prong star are refined in these two consecutive time steps.
In practice it might be more efficient to group the two regions together and treat them as one locally refined region,
thus leading to one refined discretization for the first two time steps. 
\end{remark}

\begin{table}[h]
\centering
%


\begin{tabular}{|l|c|c|c|}
\hline
&&&
\\[-1em]
 & $T_{\rm static}$&$T_{\rm Osol}$& $T_{\rm Rsol}$\\
&&&
\\[-1em]
\hline
(1) \texttt{GMRES-indy} &8.23e+2 &  3.27e+1   & 8.40e+2 \\
(2) \texttt{Direct-indy} & 8.56e+2 &  2.57e-1&  8.41e+2\\
(3) \texttt{GMRES-Local} & 8.23e+2 &   3.27e+1&  4.29e+1\\
(4) \texttt{Direct-Local} & 8.56e+2  &  2.57e-1& 1.52e+1 \\
(5) \texttt{PGMRES-Local} & 8.56e+2 &  9.40e+0&  2.43e+1\\
\hline
\end{tabular}

\caption{ Time in seconds of the construction of all necesary precomputed operators on the original discretization $T_{\rm static}$, for solving a problem that does not need local refinement $T_{\rm Osol}$ and 
for solving a problem that requires local refinement $T_{\rm Rsol}$.   Here we assume that each of the locally refined discretization is the same in size and requires the same amount of calculations to solve.}
\label{tab:star_lattice2}
\end{table}